\documentclass[preprint,authoryear,12pt]{elsarticle}

\usepackage{amssymb}
\RequirePackage[OT1]{fontenc}
\RequirePackage{amsthm,amsmath,natbib}
\RequirePackage[colorlinks,citecolor=blue,urlcolor=blue]{hyperref}
\usepackage{enumerate}

\newtheorem{Corollary}{Corollary}[section]

\newtheorem{Lemma}{Lemma}[section]
\newtheorem{theorem}{Theorem}[section]

\usepackage{amssymb,amsmath,latexsym,epsfig,hhline,graphics} 
\usepackage{amssymb}
\usepackage{graphicx}

\journal{Journal of Statistical Planning and Inference}

\begin{document}

\begin{frontmatter}



 \author{Sandra Plancade}
 \ead{sandra.c.plancade@uit.no}
 \address{Institute of Community Medicine, University of Tromso, 9037 Tromso, Norway}

\title{Adaptive estimation of the conditional cumulative distribution function from current status data}



\address{}

\begin{abstract}

Consider a positive random variable of interest $Y$ depending on a covariate $X$, and a random observation time $T$ independent of $Y$ given $X$. Assume that the only knowledge available about $Y$ is its current status at time $T$:  $\delta = 1\!\mbox{I}_{\{ Y \leq T\}}$ with $1\!\mbox{I}$ the indicator function. This paper presents a procedure to estimate the conditional cumulative distribution function $F$ of $Y$ given $X$ from an independent identically distributed sample of $(X,T,\delta)$.

A collection of finite-dimensional linear subsets of $L^2(\mathbb{R}^2)$ called models are built as tensor products of classical approximation spaces of $L^2(\mathbb{R})$. Then a collection of estimators of $F$ is constructed by minimization of a regression-type contrast on each model and a data driven procedure allows to choose an estimator among the collection. We show that the selected estimator converges as fast as the best estimator in the collection up to a multiplicative constant and is minimax over anisotropic Besov balls. Finally simulation results illustrate the performance of the estimation and underline parameters that impact the estimation accuracy.  

\end{abstract}


\end{frontmatter}

\section{Introduction}\label{section-intro}

Current status data framework refers to a situation where the observation of a positive variable of interest $Y$ called lifetime, is restricted to the knowledge of whether or not $Y$ exceeds a random measure time $T$. More precisely, one only observes the time $T$ and the ``current status'' of the system at time $T$, namely $\delta = 1\!{\mbox I}_{\{ Y \leq T \}}$ equal to 1 if $Y\leq T$ and 0 otherwise. Such data arise naturally for example in infectious disease studies, when the time $Y$ of infection is unobserved, and a test is carried out at time $T$. This framework is also named interval censoring, case 1, since the observation $(T, \delta)$ indicates whether $Y$ lies in $[0,T]$ or $(T,+\infty)$.  Contrary to other censoring framework, such as the right-censoring, the variable of interest is never directly measured. Nevertheless, as the intervals $[0,T]$ and $(T,+\infty)$ are different for each observation, quantities related to $Y$ can be inferred.

The design considered in this paper is classical in current status data framework (e.g. \cite{Cheng11}, \cite{vanderlaan06}): the survival and observation times $Y$ and $T$ are assumed to be independent given a uni-dimensional covariate $X$ directly measured. In the first place, no smoothness conditions are assumed on the variables $(X,Y,T)$.  More precisely, the following i.i.d. sample is observed
\begin{equation}\label{sample-cens-int}
(X_i, T_i, \delta _i) _{i=1, \dots ,n}
\end{equation}
with the survival times $(Y_i) _{i=1, \dots ,n}$  independent from the observation times $(T_i) _{i=1, \dots ,n}$ conditional to the covariates $(X_i) _{i=1, \dots ,n}$, and $(\delta_i ) _{i=1, \dots ,n}$ denoting the indicators $(1\!{\mbox I}_{\{ Y_i \leq T_i \}}) _{i=1, \dots ,n}$. 
We aim to estimate the cumulative distribution function (c.d.f.) $F$ of $Y$ given a covariate $X$, namely
$$F(x,y)= P[ Y \leq y | X=x].$$

Current status data have been widely studied over the past two decades. Most results about nonparametric estimation of the survival function are based on the NPMLE (Nonparametric Maximum Likelihood Estimator). \cite{groeneboom92} prove that the NPMLE is pointwise convergent  with the optimal rate $n^{-1/3}$, and \cite{vandergeer93} establishes a similar result for the $L^2 $-risk under the assumption of continuous survival time density. This unusual rate of convergence differs from the uncensored and right-censored designs, in which the distribution function can be estimated with the parametric rate of convergence $n^{-1/2}$. 

More recently, estimators developed from the NPMLE take into account the regularity of the function. \cite{hudgens07} build three estimators derived from the NPMLE, and compare their performances on simulated and real data. \cite{vanderlaan06} apply smoothing methods to the NPLME to estimate the survival function in presence of high-dimensional covariates. \cite{birge99} proposes an easily computable histogram estimator which reaches the minimax rate of convergence. The procedures developed in these papers assume that the regularity of the target function is known.

 Few papers propose estimators that automatically adapt to the unknown regularity of the function. \cite{ma} introduce a NPMLE and a least square estimator on Sobolev classes, and select the regularity parameter with a penalized criterion. \cite{comte09}  draw a parallel between current status data and bounded regression framework through a least-square estimator combined with a model selection procedure based on a more easily computable penalty function. As far as the author knows, estimation of cumulative distribution function from current status data in presence of covariates has not been handled yet in an adaptive framework. Moreover, lower-bounds on classical regularity spaces are not available. This paper aims to address these two questions. \\

 First of all, the procedure from  \cite{comte09} is extended in order to include a covariate. For every $i = 1, \dots,n$,
\begin{equation*}
\mathbb{E}[\delta _i |X_i, T_i] = \mathbb{E} \left[  1\! {\mbox I} _{\{ Y_i \leq T_i \}} |X_i, T_i \right] = \mathbb{P}[ Y_i \leq T_i |X_i, T_i ] = F(X_i,T_i).
\end{equation*}
Therefore, $F$ is the regression function of $\delta _i $ on $(X_i,T_i)$ and a least-square contrast is considered. A collection of linear spaces of $L^2(\mathbb{R}\times \mathbb{R}^+)$ is built by tensor product of classical approximation subspaces of $L^2(\mathbb{R})$ and $L^2(\mathbb{R}^+)$. Finally, a model selection estimator is computed by penalization of the least-square contrast following the model selection procedure from \cite{BM}.

The loss function usually considered in model selection procedures is associated to the $L^2$-norm. Nevertheless, the least-square contrast which quantifies the goodness of fit on the observations provides a natural control of the risk associated to the empirical norm $\|t \|_n ^2 = (1/n) \sum _{i=1}^n t^2(X_i,T_i)$ which represents an intermediate result  in the computation of the integrated risk in most model selection procedures (see e.g. \cite{massart}, Comments of Theorem 8.3). In this paper, we have chosen to first display the control of the empirical risk conditional to the observations $\{ (X_i, T_i)\}$ and then derive the result for the $L^2$-risk under additional assumptions in order to emphasize the role played by each hypothesis. Besides, the fixed design context with non-random observations $\{ (X_i, T_i)\}$ can be directly handled.

 On complementary interest to the adaptivity on decomposition bases, is the study of the rate of convergence on spaces of regularity. Anisotropic Besov spaces, which account  for different regularities of the target function $F$ with respect to the covariate $X$ and the survival time $Y$, are considered. 
On the one hand, the rate of convergence of our adaptive estimator is evaluated, exhibiting classical rates in a regression context, but unusual in c.d.f. estimation. On the other hand, the minimax rate of convergence is considered as a criterion to evaluate the fundamental limit of the current status framework. The combination of these two analyses assesses the optimality of our estimator.  
 Results in a context without covariate can be derived as a particular case, proving the optimality of the estimators proposed by \cite{ma} and \cite{comte09}. 

 Furthermore, we analyze the performance of our estimator on simulated data. In a right-censoring context, the rate of censoring defined as the proportion of censored observations is a well-identified parameter of the quality of the estimation. Since the survival time is never observed in current status data, this quantity makes no sense. Nevertheless, in our numerical study, we exhibit a parameter which impacts the estimation from current status data, and may represent the counterpart of the rate of censoring on right-censored data.  \\

The paper is organized as follows. Section \ref{sec-notations} introduces the notations. The estimation procedure is presented in Section \ref{sec-estim}. The upper-bounds for the empirical and integrated loss functions are presented in Section \ref{sec-upper}. The minimax study over anisotropic Besov balls is conducted in Section \ref{sec-minimax}. Section \ref{sec-simus} displays the performances of the estimator on simulated data. General comments and perspectives for further developments are presented in Section \ref{sec-comments}. The proofs are gathered in Section \ref{sec-proofs}.

\section{Notations} \label{sec-notations}

\subsection{Notations}

Let {\bf X} and {\bf T} denote the vectors $(X_1,\dots,X_n)$ and $(T_1,\dots,T_n)$, and for every $t \in L^2(\mathbb{R}\times \mathbb{R}^+)$, let $t({\bf X, T}) = (t(X_1,T_1), \dots, t(X_n,T_n) )$. 
For every  $(V_i , W_i)_{i=1\dots, n}$ i.i.d. random variables, $f_V $ denotes the density of $V_i$ and $f_{(V,W)} $ the density of the couple $(V_i,W_i)$.

 Let $A=A_1 \times A_2$ be a subset of $\mathbb{R}\times \mathbb{R}^+$ with $A_1$ and $A_2$ compact intervals. For every $t$, $s\in L^2 (A)$, we  define the empirical norm and scalar product relative to the sample $(X_i,T_i)_{i=1,\dots,n}$:

$$\langle s, t \rangle _n = \frac{1}{n} \sum _{i=1}^n s(X_i, T_i)t(X_i,T_i) \quad \text{and}  \quad \| s\| _n ^2 = \langle s,s \rangle _n $$
and the associated integrated norm and scalar product:
$$\langle s,t \rangle _{f_{(X,T)}} = \int _{x\in A_1} \int _{u\in A_2} s(x,u) t(x,u) f_{(X,T)} (x,u) dx du \quad \text{and}  \quad \| s\|  _{f_{(X,T)}}^2 = \langle s,s \rangle  _{f_{(X,T)}} .$$
Finally, for every finite set $\mathcal{S}$, $| \mathcal{S}|$ denotes the cardinal of  $\mathcal{S}$.

\section{Model selection estimator}\label{sec-estim}

\subsection{Collection of linear models}

We build a collection of linear subspaces of $L^2(A)$, referred to as \textit{models}, by tensor product of finite-dimensional linear subspaces of $L^2(A_1)$ and $L^2(A_2)$. For $j=1$ or $2$, let 
\begin{equation}\label{coll-Mn} 
\mathcal{M}_n ^{(j)} = \{ S_{m_j}^{(j)} , m_j\in I_n^{(j)} \} 
\end{equation}
be a collection of linear subspaces of $L^2(A_j)$ with finite dimension $D_{m_j}^{(j)}$. Let  $I_n =I_n^{(1)}\times I_n^{(2)}$. For every $m= (m_1 , m_2 ) \in I_n$, we define

\begin{equation}\label{tensor}
S_m = S_{m_1}^{(1)} \otimes S_{m_2}^{(2)} =  \left\{ t :A \rightarrow \mathbb{R},\quad  t(x,y)= \sum _{(k,l)\in J_m}a_{k,l} \phi _k^{m_1}(x)  \psi _l^{m_2}(y) \right\}
\end{equation}
with $(\phi _k ^{m_1}) _{k=1, \dots , D_{m_1}^{(1)}}$ an orthonormal basis of $S_{m_1}^{(1)}$, $(\psi _l ^{m_2}) _{l=1, \dots , D_{m_2}^{(2)}}$ an orthonormal basis of $S_{m_2}^{(2)}$ and
\begin{equation}\label{Jm}
J_m= \left( (1,1) , \dots , (1, D_{m_2}^{(2)}), (2,1) , \dots , (2, D_{m_2}^{(2)}),  \dots, (D_{m_1}^{(1)},1), \dots, (D_{m_1}^{(1)}, D_{m_2}^{(2)}) \right)
\end{equation}
For every $m= (m_1, m_2)$, $S_m$ is a linear subspace of $L^2(A)$ with dimension $D_m = D_{m_1}^{(1)} D_{m_2}^{(2)}$. Then the collection of models is defined as $\mathcal{M}_n = \{ S_m, m\in I_n \} .$
 Once fixed the collection $\mathcal{M}_n$, a bijection is established between the indexes $\{m \in I_n\}$ and the linear spaces $\{S_m \in \mathcal{M}_n\}$. Therefore, the term "model" will indistinctly refer either to $m$ or $S_m$. 
We assume that $\mathcal{M}_n^{(1)}$ and $\mathcal{M}_n^{(2)}$ satisfy the following assumption.\\

${\bf (H) } $ Let $j=1$ or $2$. For every $b>0$, there exists a constant $B_j$ such that
\begin{equation*}
\sum _{m_j\in I_n^{(j)}} \exp \left(- b \sqrt{D^{(j)}_{m_j}}\right) \leq B_j, \quad \forall n\in \mathbb{N}^*.
\end{equation*}

\subsection{Least-square estimators}\label{sec-leastsquare}
For every $t \in L^2(A)$, we consider the least-square contrast

\begin{equation}\label{gamman}
\gamma _n (t) = \frac{1}{n} \sum _{i=1}^n \left( \delta _i - t(X_i, T_i) \right)^2.
\end{equation}
For every model $S_m \in \mathcal{M}_n$, we denote by $\widehat{F}_m$ the minimum contrast estimator on $S_m$:
\begin{equation}\label{Fhatm}
\widehat{F}_m= \arg\min _{t\in S_m } \gamma _n (t) = \sum _{(k,l)\in J_m} \widehat{a}_{k,l} \phi _k^{m_1} (.) \psi _l ^{m_2} (.)
\end{equation}
Equation (\ref{Fhatm}) amounts to state that the partial derivatives $\partial \gamma _n (\widehat{F}_m)/ \partial a _{k,l} $ are equal to zero for every $(k,l) \in J_m$. Therefore, the  column vector of coefficients $\widehat{A}_m=[ \widehat{a}_{k,l}]_{(k,l)\in J_m}$ satisfies

\begin{equation}\label{eq-Am}
\widehat{G}_m \widehat{A}_m=\widehat{V}_m
\end{equation}
with
$$\widehat{G}_m= \left[ \frac{1}{n} \sum _{i=1}^n \phi _k ^{m_1} (X_i) \psi _l^{m_2} (T_i)\phi _{k'} ^{m_1} (X_i) \psi _{l'}^{m_2} (T_i) \right] _{((k,l), (k',l'))\in J_m^2 }$$
the $D_m \times D_m$-square  Gram matrix related to $\{ \phi _k^{m_1} \psi _l^{m_2} \} _{(k,l)\in J_m}$ for the scalar product $\langle.,.\rangle _n$ and
$$\widehat{V}_m=\left[ \frac{1}{n} \sum _{i=1}^n \delta _i \phi _k^{m_1} (X_i) \psi _l ^{m_2} (T_i)  \right] _{(k,l) \in J_m} $$
 a column vector. \\

\textbf{Comment on existence and unicity of $\widehat{F}_m$.} 
Equation (\ref{eq-Am}) has a unique solution $\widehat{A}_m$ if and only if the matrix $\widehat{G}_m$ is invertible. Nevertheless, we show that the estimator $\widehat{F}_m$ exists almost surely and is uniquely defined on the set $({\bf X,T})$. 
 
 Consider the observed sample (\ref{sample-cens-int}). Let $\mathcal{S}$ be a finite dimensional linear subset of $L^2(A)$ and $\Phi$ be the linear application from $\mathcal{S}$ to $\mathbb{R}^n$ defined as $\Phi(t) = t({\bf X,T})$. We denote by $\mathcal{S}^*$ the linear subset of $\mathcal{S}$ such that:
 
 \begin{equation}\label{eq-Sstar}
  \mathcal{S}^* \oplus \text{Ker} (\Phi) =\mathcal{S}
  \end{equation}
 with $\text{Ker} (\Phi)$ the null space of $\Phi$. $\mathcal{S}^*$ is the largest linear subset of $\mathcal{S}$ on which $\|.\|_n$ satisfies the separation property: $\|t\|_n =0 \Rightarrow t=0$.  The restriction $\Phi^*$ of $\Phi$ to $\mathcal{S}^*$ draws a bijection from $\mathcal{S}^*$ into the image $\Phi(\mathcal{S})$ of $\Phi$. Let $Z_{\mathcal{S}}$ be the projection of $(\delta _1,\dots,\delta _n )$ on $\Phi(\mathcal{S})$. There exists a unique function $\widehat{F}^*_{\mathcal{S}} \in \mathcal{S}^*$ such that $\Phi(\widehat{F}^*_{\mathcal{S}}) =  Z_{\mathcal{S}}$. Furthermore, $\widehat{F}^*_{\mathcal{S}} = \arg \min _{t \in \mathcal{S}} \gamma_n(t) $ and for any minimiser $\widehat{F}_{\mathcal{S}}$ of $\gamma _n $ over $\mathcal{S}$, $\widehat{F}^*_{\mathcal{S}}$ and $\widehat{F}_{\mathcal{S}}$ are equal on $({\bf X,T})$.

\subsection{Model selection estimator}

For every  $m \in \mathcal{M}_n$, let $F_m$ be a minimizer of $\| F-t \| _n^2$ for $t \in S_m$. The empirical risk of the least-square estimator $\widehat{F}_m$ has the following bias-variance decomposition:

$$\mathbb{E}\left[ \| \widehat{F}_m - F \|_n^2| ({\bf X,T}) \right] = \mathbb{E}\left[ \| \widehat{F}_m - F_m \|_n^2| ({\bf X,T}) \right] + \| F_m -F \| _n^2.$$
The bias term 
 $\| F_m -F \| _n^2 $  is estimated by $\gamma _n (\widehat{F}_m)$.
 
 Let $S_m^*$ be defined by (\ref{eq-Sstar}), $\widehat{F}_m^*$ and $F_m^*$ be respectively the minimisers of $\gamma _n$ and $\|.-F\|_n$  over $S_m^*$. Then $\widehat{F}_m^*({\bf X,T})$ and $F_m^*({\bf X,T})$ are respectively the projection of $(\delta _1,\dots,\delta _n )$ and $F({\bf X,T})$ on $S_m^*$ for the norm $\|.\|_n$. Thus, let $\left\{ \varphi _{\lambda}, \lambda \in J_m^* \right\}$ be a $\|.\|_n$-orthonormal basis of $S_m^*$,

 $$\| \widehat{F}_m^*- F_m^*\|_n^2 = \| \widehat{F}_m- F_m\|_n^2 = \sum _{\lambda \in J_m^*} \left( \nu _n (\varphi_{\lambda}) \right)^2$$
with
\begin{equation}\label{eq-nu}
\nu _n(t) = \frac{1}{n} \sum _{i=1}^n ( \delta _i - F(X_i,T_i))t(X_i,T_i).
\end{equation}
Moreover, the $\{ \delta _i \} _{i=1,\dots,n}$ are independent conditional to $({\bf X,T})$ and $ \mathbb{E} [(\delta _i - F(X_i,T_i))^2 |{\bf X,T}] \leq 1/4$ as variance of a Bernouilli variable, which entails:
\begin{equation}\label{eq-bias}
\mathbb{E} \left[ (\nu _n (t))^2  |{\bf X,T} \right] \leq \frac{1}{4n^2} \sum _{i=1}^n t^2(X_i, T_i) = \frac{\| t\| _n^2}{4n} , \quad \forall t \in L^2(A).
\end{equation}
 Therefore, the variance term $\mathbb{E}\left[ \| \widehat{F}_m - F_m \|_n^2| ({\bf X,T}) \right]$ is upper bounded by $(1/4) D_m/n$ and a model $\widehat{m}$ which minimizes the estimated bias-variance sum is selected:
 $$\widehat{m} = \arg \min _{m \in I_m} \left[ \gamma _n(\widehat{F}_m ) + pen(m) \right]$$
 with $pen(m) = (\theta/4) D_m/n$ for some numerical constant $\theta >1$. \\

Besides, the cumulative distribution function $F$ lies in $[0,1]$. We use this information to improve the estimation by constraining the values of our estimator to remain in this interval. More precisely we consider the estimator $\tilde{F}_{\widehat{m}}$ with

\begin{equation}\label{Ftilde-def}
\begin{array}{l}
\\
\tilde{F}_m (x,u) = \\
\\
\end{array}
\left\{ \begin{array}{l}
0 \quad \text{if} \quad \widehat{F}_m (x,u) <0 \\
1   \quad \text{if} \quad \widehat{F}_m (x,u) >1 \\
\widehat{F}_m(x,u) \quad \text{otherwise}.
\end{array}\right.
\end{equation}

\noindent \textbf{Comments}

\begin{enumerate}[1)]

\item The restriction imposed by (\ref{Ftilde-def}) 
brings an improvement in terms of risk. Indeed, for every $(x,u) \in A$, 
$$|\tilde{F} _m (x,u) - F(x,u) | \leq | \widehat{F}_m(x,u) - F(x,u) |, \quad \forall m \in \mathcal{M}_n $$
almost surely. In particular, $\|\tilde{F}_m - F\|_n^2 \leq \|\widehat{F}_m - F\|_n^2 $.
Thus, any upper-bound of $\mathbb{E}[ \| \widehat{F}_m -F \| _n ^2 |  {\bf X,T}  ]$ is also an upper-bound of  $\mathbb{E}[ \| \tilde{F}_m -F \| _n ^2 | {\bf X,T} ]$. \\

\item The convergence results presented in this paper are valid for any $\theta>1$, but in practical implementation a value  has to be fixed. This issue arises in most model selection procedure. It can be either calibrated on simulated data from a large number of examples,  determined through a cross-validation procedure, or chosen a priori. In particular,  \cite{massart-pen-ct} proposes an heuristic which validates the value $\theta = 2$ in a general context. \\

\item As $F$ is bounded by $1$, the constant involved in the penalty is purely numerical contrary to many frameworks where the penalty includes unknown parameters which have to be estimated. 

\end{enumerate}

\section{Non-asymptotic upper-bounds}\label{sec-upper}

\subsection{Oracle inequality for the empirical risk}

The minimization of the least-square contrast (\ref{gamman}) determines the function which offers the best fit on the sample $({\bf X,T}) $. Therefore,  the empirical loss function $\| \tilde{F}_{\widehat{m}} -F \| _n ^2$ which measures the error of estimation on the set $({\bf X,T}) $ appears as the natural risk for the  estimator $\tilde{F}_{\widehat{m}}$.

\begin{theorem}\label{T1}
Assume that ${\bf (H)}$ holds.
There exist numerical constants $C_1$ and $C_2$ such that

\begin{equation}\label{eq-T1-chap5}
\mathbb{E}\left[ \| \tilde{F}_{\widehat{m}} -F \| _n ^2 |  {\bf X,T} \right]  \leq C_1 \inf _{m\in I_n} \left\{  \inf _{t\in S_m } \| F-t \| _n ^2 + pen (m) \right\} + \frac{C_2}{n} \quad a.s.
\end{equation}

\end{theorem}

 \noindent \textbf{Comments on Theorem \ref{T1}}

\begin{enumerate}[1)]

\item  Theorem $\ref{T1}$ indicates that the model selection estimator  realizes the trade-off between the bias and the penalty, up to a multiplicative constant. Nevertheless, since the penalty has only been stated as an upper-bound of the variance, we did not strictly proved that the sum $(  \inf _{t\in S_m }\| F-t \| _n ^2 + pen (m))$ has the order of the bias-variance sum. The minimax study conducted in Section \ref{sec-minimax} is necessary to ensure that the penalty has the order of the variance, and therefore assess the adaptivity of the estimator.  \\

\item The result of Theorem $\ref{T1}$ holds for any value of the sample $({\bf X,T}) $ , regardless of its distribution. In particular, it handles the situation where the observation times ${\bf T}$ and the covariates ${\bf X}$ are fixed by the experimental design. \\

\end{enumerate}
 
\subsection{Oracle inequality for the $L^2$-risk}\label{section-int-oracle}

Theorem \ref{T1} states the adaptivity of the model selection estimator for the empirical risk which naturally arises from the least-square contrast. Nevertheless, the loss function considered in Theorem \ref{T1} is data dependent, and  the control of the classical $L^2$-risk provides more interpretable results. In this section, we derive an oracle inequality for the $L^2$-risk  under the following assumptions. \\

${\bf (A_1)}$ For $i=1,\dots, n$, $(X_i,T_i)$ has a density $ f_{(X,T)}$ bounded by  $h_0>0$ and $h_1 <+\infty$ on set $A$:
$$ h_0 \leq f_{(X,T)}(x,u) \leq h_1, \quad \forall (x,u) \in A.$$

${\bf (A_2)} $  For $j=1$ and 2, 
\begin{equation}\label{card-poly}
\left| \mathcal{M}_n^{(j)} \right| \leq P^{(j)}(n), \quad \forall n \in \mathbb{N} 
\end{equation}
for some polynomial $P^{(j)}$. Moreover, there  exists a model $S_n^{(j)} \in \mathcal{M}_n^{(j)}$ of dimension $N_n^{(j)}$
such that, for every $m_j\in I_n^{(j)} $, $S_{m_j}^{(j)} \subset S_n^{(j)}$ and $N_n^{(1)} N_n^{(2)} \leq \sqrt{n}/ \log n.$ \\

${\bf (A_3)}$ There exists a positive constant $K_1$ (resp. $K_2$) such that, for every $m_1 \in I_n^{(1)}$ (resp. $m_2 \in I_n^{(2)}$),
$$\sup _{x\in A_1}  \sum _{k=1}^{D_{m_1}^{(1)}} (\phi _k^{m_1} (x))^2  \leq K_1 D_{m_1}^{(1)} \quad \left( \text{resp.     } \sup _{u\in A_2}  \sum _{l=1}^{D_{m_2}^{(2)}} (\psi _l^{m_2} (u))^2  \leq K_2 D_{m_2}^{(2)} \right) .$$

\medskip

Assumption ${\bf (A_2)}$ refers to the number of models in the collection, whereas ${\bf (A_3)}$ is related to the nature of the models. \\

\textbf{Examples of models.} Assumption ${\bf (A_3)}$ holds when the collections $\mathcal{M}_n^{(1)} $ and $\mathcal{M}_n^{(2)} $  consist of the following models. 

\begin{enumerate}[(i)]

\item $S_{m_j}^{(j)}$ is the set of piecewise polynomials with maximum degree $s_j$ and step $\text{lg}(A_j)/ D_{m_j}^{(j)}$ where $\text{lg}(A_j)$ denotes the length of $A_j$.\\

\item $S_{m_j}^{(j)} =\text{vect}\{ \chi _{l,k} , l\leq m_j , k\in \mathbb{Z} \}$ where $\chi$ is a bounded mother wavelet with regularity $s_j $, $\chi _{l,k} (x) = 2^{l/2} \chi (2^l x-k) $ and $D_{m_j}^{(j)}= 2^{m_j}$.\\

\item $S_{m_j}^{(j)}$ is the set of trigonometric polynomials with maximum degree $D_{m_j}^{(j)}$.\\

\end{enumerate}

 Theorem $\ref{T1}$ leads to the following result.  

\begin{theorem}\label{T2-cens-int}
Assume that ${\bf (H)}$, ${\bf (A_1)}$, ${\bf (A_2)}$ and ${\bf (A_3)}$ hold. Then
\begin{equation}\label{eq-C1}
\mathbb{E} \left[ \| \tilde{F}_{\widehat{m}} -F \| ^2 \right] \leq C_3 \inf _{m\in I_n} \left\{  \inf _{t\in S_m } \| F-t \|  ^2 + pen (m) \right\} + \frac{C_4}{n}
\end{equation}
for some constant $(C_3,C_4)$ depending on $(h_0,h_1,K)$.

\end{theorem}

\section{Minimax study  on Besov balls}\label{sec-minimax}
Theorem \ref{T2-cens-int} states the adaptivity of the model selection estimator i.e. its optimality among the collection of estimators. The minimax rate of convergence, defined as the rate of the best estimator for the worst function among a given class of regularity, assesses optimality in a more general sense. In this section, we prove that  the model selection estimator is minimax over anisotropic Besov balls.

\subsection{Minimax rate of convergence}

Let $\mathcal{F}$ be a set of conditional cumulative distribution functions on $A$ and $(r_n )_{n\in \mathbb{N}}$ a sequence  of positive numbers. The sequence $(r_n )_{n\in \mathbb{N}}$ is the minimax rate of convergence of $F$ over $\mathcal{F}$, defined up to a constant,  if there exist two constants $c$ and $C$ such that
$$c \leq \inf _{\widehat{F}_n} \sup _{F\in \mathcal{F}} \left( r_n ^{-1}\mathbb{E} [ \| \widehat{F}_n - F \|^2 ] \right)\leq C$$
with the infimum taken over all possible estimators $\widehat{F}_n$. The minimax rate corresponds to the rate of convergence of the best possible estimator of the function $F$, based on the sample $(X_i,T_i,\delta_i)_{i=1,\dots,n}$. 

\subsection{Definition of anisotropic Besov balls}

We recall the definition of two-dimensional anisotropic Besov spaces stated for example in \cite{hochmuth}.
Let $f\in L^2 (A)$. For $j=1$ or $2$, $r\in \mathbb{N}^*$ and $h>0$, let 
$$\Delta _{h,j}^r (f)(x,y)= \sum _{k=0}^r \left( _k ^r \right) (-1)^{r-k} f((x,y)+khe_j)$$
be the directional partial difference operator for every $(x, y) \in A _{h,j}^r $ where
$A _{h,j}^r = \left\{ (x,y)\in A, (x,y)+ rhe_j \in A \right\}$
and $(e_1,e_2)$ is the canonical basis of $\mathbb{R}^2$. For $t>0$, let
$A _{r,j} (f,t, A) = \sup _{|h| \leq t } \| \Delta _{h,j}^r (f)(x,y) \| _{L^2(A _{h,j}^r)}$
be the directional modulus of smoothness for the $L^2$-norm.
Let $\beta =(\beta _1, \beta _2) \in (\mathbb{R}_+^*) ^2 $ and $r_j= \lfloor \beta _j \rfloor +1$.  We define the anisotropic Besov space of parameters $(\beta, 2, \infty)$ as

$$\mathcal{B}^{\beta}_{2, \infty}(A)  = \left\{ f \in L^2 (A) , | f | _{ \mathcal{B}^{\beta}_{2, \infty}(A)} <+\infty \right\}$$
with
$$| f | _{ \mathcal{B}^{\beta}_{2, \infty}(A)} = \sup _{t>0} \left[ t^{-\beta _1} A _{r_1,1} (f,t,A) + t^{-\beta _2}  A _{r_2,2} (f,t,A)\right] .$$
The Besov norm on $\mathcal{B}^{\beta}_{2, \infty}(A)$ is defined by
$\| f \| _{ \mathcal{B}^{\beta}_{2, \infty}(A)}=| f | _{ \mathcal{B}^{\beta}_{2, \infty}(A)} + \| f\|$
and for $L>0$, we consider the anisotropic Besov ball:
$$\mathcal{B}^{\beta}_{2, \infty} (A, L) = \left\{ f \in \mathcal{B}^{\beta}_{2, \infty}(A), \| f \| _{ \mathcal{B}^{\beta}_{2, \infty}(A)}  \leq L  \right\}.$$

\subsection{Upper bound of the $L^2$-risk on anisotropic Besov balls} The following corollary of Theorem \ref{T2-cens-int} assesses the rate of convergence of the model selection estimator on anisotropic Besov balls. 

\begin{Corollary}\label{C2-cens-int}

Assume that $F\in \mathcal{B}^{\beta}_{2, \infty} (A, L)$ for some $L>0$ and $\beta = (\beta _1, \beta _2 ) \in ( \mathbb{R}_+^*)^2$  and  $\mathcal{M}_n$ is set up from models $(i)$, $(ii)$ with $s_j > \beta _j -1$, or $(iii)$. 
Then
\begin{equation}\label{L-bias}
\inf _{t\in S_m} \| F-t\| \leq C_0 \left( (D_{m_1}^{(1)})^{-\beta _1} + (D_{m_2}^{(2)})^{-\beta _2}\right).
\end{equation}
for some  positive constant $C_0$. Moreover assume that ${\bf (A_1)}$ and ${\bf (A_2)}$ hold, and
\begin{equation}\label{Nnj-chap5}
N_n^{(j)} \leq \left( \frac{n}{\log ^2 n} \right)^{1/4} \quad \text{for} \ \ \ j=1,2.
\end{equation}
Then

$$\mathbb{E} \left[ \| \tilde{F}_{\widehat{m}} - F \|  ^2 \right] \leq C _5 n ^{- \overline{\beta} /(\overline{\beta} +1)}$$
for some positive constant $C _5$, with $\overline{\beta} = 2\beta _1 \beta _2 /(\beta _1 + \beta _2)$ the harmonic mean of $(\beta _1, \beta _2)$.

\end{Corollary}

The upper-bound of the  bias term (\ref{L-bias}) is computed in \cite{claire07}  
 based on results from \cite{hochmuth} and \cite{nikolskii}. 
Plugging this result into Theorem \ref{T2-cens-int} provides the rate of convergence of the bias-variance sum:
\begin{equation}\label{eq-trade}
   \inf _{t\in S_m } \| F-t \|  ^2 + pen (m)  \leq  C_6
\left\{ \left((D_{m_1}^{(1)})^{-2\beta _1} + (D_{m_2}^{(2)})^{-2\beta _2}\right) + \frac{ D_{m_1}^{(1)} D_{m_2}^{(2)}}{n} \right\} 
\end{equation}
and the trade-off between the two terms in the right hand of (\ref{eq-trade})  is achieved for a model $(m_1,m_2)$ satisfying:
$$ D_{m_1}^{(1)} \propto n^{\beta _2/(\beta _1 + \beta _2 + 2 \beta _1 \beta _2)} \quad \text{and} \quad 
D_{m_2}^{(2)} \propto n^{\beta _1/(\beta _1 + \beta _2 + 2 \beta _1 \beta _2)},$$
with $\propto$ denoting proportionality. Moreover, Assumption (\ref{Nnj-chap5}) guarantees the existence of a model $(\overline{m}_1,\overline{m}_2)$ such that
$$ 1 \leq D_{\overline{m}_1}^{(1)} n^{-\beta _2/(\beta _1 + \beta _2 + 2 \beta _1 \beta _2)} \leq 2 \quad \text{and} \quad 1 \leq D_{\overline{m}_2}^{(2)} n^{-\beta _1/(\beta _1 + \beta _2 + 2 \beta _1 \beta _2)} \leq 2$$
for   $\beta _1, \ \beta _2 >1$. Therefore,
 $$ \left( (D_{\overline{m}_1}^{(1)})^{-2\beta _1} + D_{\overline{m}_2} ^{(2)})^{-2\beta _2}\right) +\frac{ D_{\overline{m}_1}^{(1)} D_{\overline{m}_2}^{(2)}}{n}\propto C_7 n ^{- \overline{\beta} /(\overline{\beta} +1)},$$
which concludes the proof of Corollary \ref{C2-cens-int}. $\Box$\\

\noindent \textbf{Comments.} Condition $\beta _1$, $\beta _2 >1$ in Corollary $\ref{C2-cens-int}$ can be extended to
$$(\beta _1, \beta _2) \in ( \beta _1^* , + \infty) \times (\beta _2^* , + \infty )$$ for a known couple $(\beta _1^*, \beta _2^*)$ which satisfies $\overline{\beta}^* \geq 1$, with $\overline{\beta }^*$ the harmonic mean of $\beta _1^*$ and $\beta _2^*$, by considering $N_n^{(1)}$ and $N_n^{(2)}$ such that
\begin{equation*}
\begin{array}{l}
N_n^{(1)} \leq ( \log n)^{-1/2} n^{\beta _1^*/(\beta _1^* + \beta _2^* + 2 \beta _1^* \beta _2^*)} \quad \text{and} \quad N_n^{(2)} \leq ( \log n)^{-1/2}n^{\beta _2^*/(\beta _1^* + \beta _2^* + 2 \beta _1^* \beta _2^*)}.
\end{array}
\end{equation*}
This alternative assumption takes into account a priori knowledge on the regularity of $F$, through an appropriate choice of $(N_n^{(1)}, N_n^{(2)})$. Thus, the estimation would be optimized by considering a smaller maximum size of models $(N_n^{(j)})$ in the direction where $F$ is more regular. 

\subsection{Lower bound of the minimax rate of convergence}

\begin{theorem}\label{T-low}
Let $\beta =(\beta _1, \beta _2) \in (0,+\infty) \times (1, + \infty)$.
Assume that $h_1 = \| f_{(X,T) } \|_{\infty} <+\infty$, then there exists a constant $c$ which depends on $(\beta, L, h_1 ) $ such that
$$\inf_{\widehat{F}_n } \sup _{F\in \mathcal{B}^{\beta}_{2, \infty} (A , L) } \mathbb{E} \left[ n^{-\overline{\beta}/(\overline{\beta} + 1)} \| \widehat{F}_n - F \| ^2 \right] \geq c.$$
with the infimum taken over all possible estimators $\widehat{F}_n$ based on sample (\ref{sample-cens-int})
\end{theorem}

\subsection{Optimality of the model selection estimator in the minimax sense}

Let  $\beta _1$, $\beta _2 >1$. By combining the results from Theorem \ref{T2-cens-int} and \ref{T-low}, we obtain:

\begin{equation*}
\hspace{-0.6cm}c\leq \inf_{\widehat{F}_n }\sup _{F\in \mathcal{B}^{\beta}_{2, \infty} (A , L) } \mathbb{E} \left[ n^{-\overline{\beta}/(\overline{\beta} + 1)} \| \widehat{F}_n - F \| ^2 \right] \leq \sup _{F\in \mathcal{B}^{\beta}_{2, \infty} (A , L) } \mathbb{E} \left[ n^{-\overline{\beta}/(\overline{\beta} + 1)} \| \widehat{F}_{\widehat{m}} - F \| ^2 \right]  \leq C
\end{equation*}
for some positive constants $c$ and $C$. Therefore,  the minimax rate of convergence over $\mathcal{B}^{\beta}_{2, \infty} (A , L)$ is $n^{-\overline{\beta}/(\overline{\beta} +1)}$, and the model selection estimator $\tilde{F}_{\widehat{m}}$ is optimal in the minimax sense. \\

\noindent \textbf{Comments.} 

\begin{enumerate}
\item In the context where $Y$ and $T$ do not depend on a covariate, the minimax rate of convergence on the Besov balls of regularity $\beta $ is $n^{-2\beta/(2\beta+1)}$. Heuristically, this situation corresponds to an "infinite" regularity $\beta _1$ with respect to the first variable, the harmonic mean $\overline{\beta}$ being equal to $2\beta_2$. The formal proof for the lower bound is given in Section \ref{sec-proofs}, and the upper bound is provided by the estimator from \cite{comte09}. 

\item The rate $n^{-\overline{\beta}/(\overline{\beta} + 1)}$ corresponds to the minimax rate in multivariate regression estimation (e.g. \cite{Gaiffas11}) and differs from the rate $n^{-2\beta_1/(2\beta_1+1)}$ obtained on right-censored data, with $\beta _1$ the c.d.f. with respect to the covariate $X$ (\cite{Brunel07}). This result is coherent with the results from \cite{comte09} for independent survival and observation times: whereas the c.d.f. is estimated at the nonparametric rate $1/n$ from right-censored data, in the interval censoring framework the minimax rate over Besov balls is the classical regression estimation rate $n^{2\alpha/(2\alpha+1)}$ with $\alpha$ the regularity of the c.d.f. 

\end{enumerate}
\section{Numerical study}\label{sec-simus}

\subsection{Numerical implementation}\label{section-compute-details}

  We  implement the model selection estimator of the conditional c.d.f. on a collection of models generated by histogram bases. 
The constant $\theta$ in the penalty is set to $2$ refering to \cite{massart-pen-ct}. Moreover, our studies indicate that the results are robust with respect to the value of the constant (not shown).

\subsection{Two examples: multiplicative and additive effects of the variable $X$}\label{sec-mod12}
\mbox{}

\vspace{0.2cm}
\noindent \textbf{Model 1: additive effect of $X$.} We consider the following distribution of $(X,Y,T)$:

$$\begin{array}{l}
\\
\text{(Mod 1)}\\
\\
\end{array}
\left\{\begin{array}{l}
X\sim \Gamma (k=1,\theta=1) \\
Y = X+ \varepsilon  \quad \text{with} \quad  \varepsilon \sim \Gamma (k=3,\theta=2) \\
T = X+ \varepsilon'  \quad \text{with} \quad   \varepsilon ' \sim \Gamma (k=3,\theta=2) \\
\end{array} \right.$$
where $\Gamma (k,\theta)$ denotes the gamma distribution with shape $k$ and scale $\theta$. The conditional cumulative distribution function $F(x,y) = \phi(y-x, k=3,\theta=2)$, where $\phi$ denotes the density of the gamma distribution, is estimated on the set $A=[0.1,3]\times [3.1,19]$. The density $f_{(X,T)}$ is lower bounded by $h_0=3.10^{-5}$ on A. \\

\noindent \textbf{Model 2: multiplicative effect of $X$.} We consider the following distribution of $(X,Y,T)$:

$$\begin{array}{l}
\\
\text{(Mod 2)}\\
\\
\end{array}\left\{\begin{array}{l}
X\sim \Gamma (k=1.5,\theta=2) \\
Y = X \times \varepsilon \quad \text{with} \quad  \varepsilon \sim \Gamma (k=3,\theta=1) \\
T = X\times \varepsilon' \quad \text{with} \quad  \varepsilon ' \sim \Gamma (k=3,\theta=1) \\
\end{array}
\right.$$
 The conditional cumulative distribution function $F(x,y) = \phi(y/x, k=3,\theta=1)$ is estimated on the set $A=[1,10]\times [1,20]$. The density $f_{(X,T)}$ is lower bounded by $h_0=10^{-7}$ on A. \\

\begin{figure}
\begin{tabular}{cccc}
\hspace{-1.5cm}\includegraphics[width=4cm,trim= 20mm 20mm 10mm 20mm]{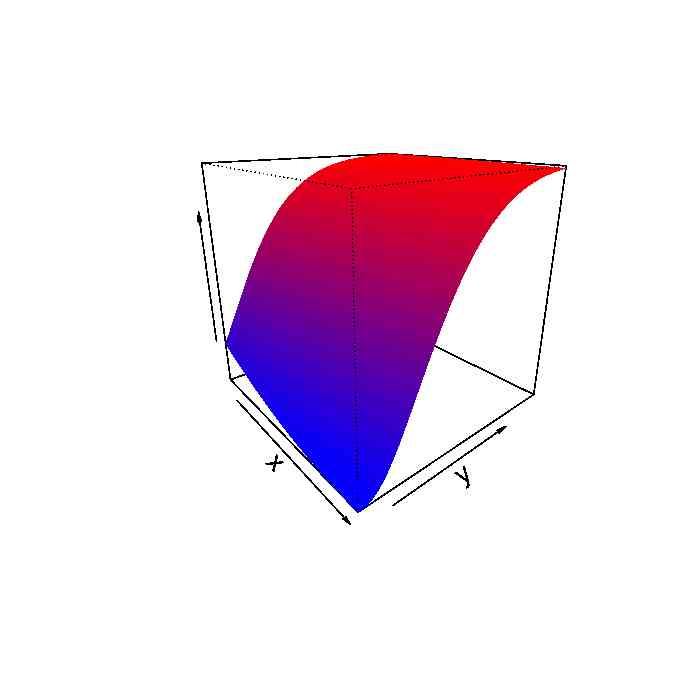} &
\includegraphics[width=4cm,trim= 20mm 20mm 10mm 20mm]{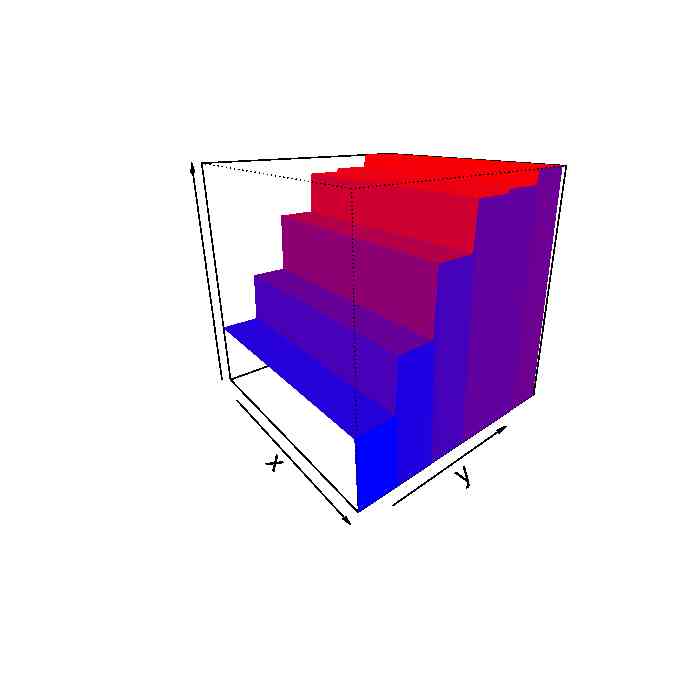} &
 \includegraphics[width=4cm,trim= 20mm 20mm 10mm 20mm]{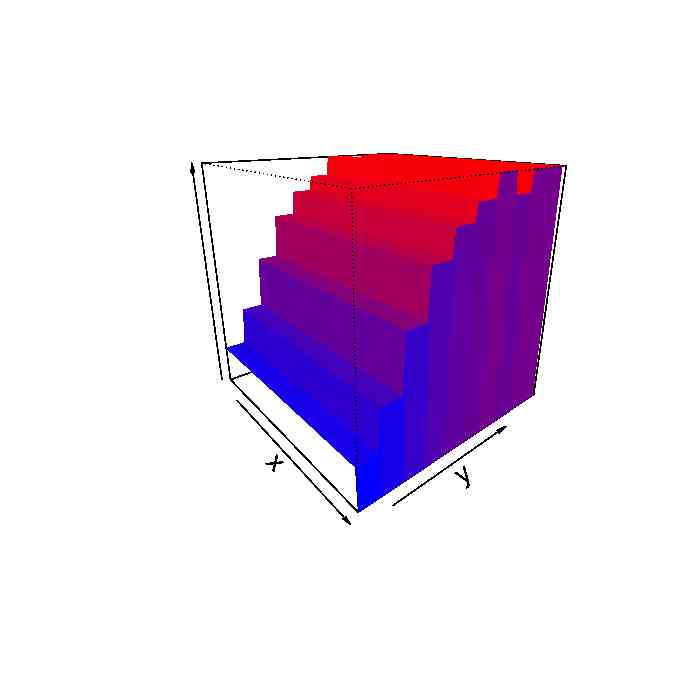} &
 \includegraphics[width=4cm,trim= 20mm 20mm 10mm 20mm]{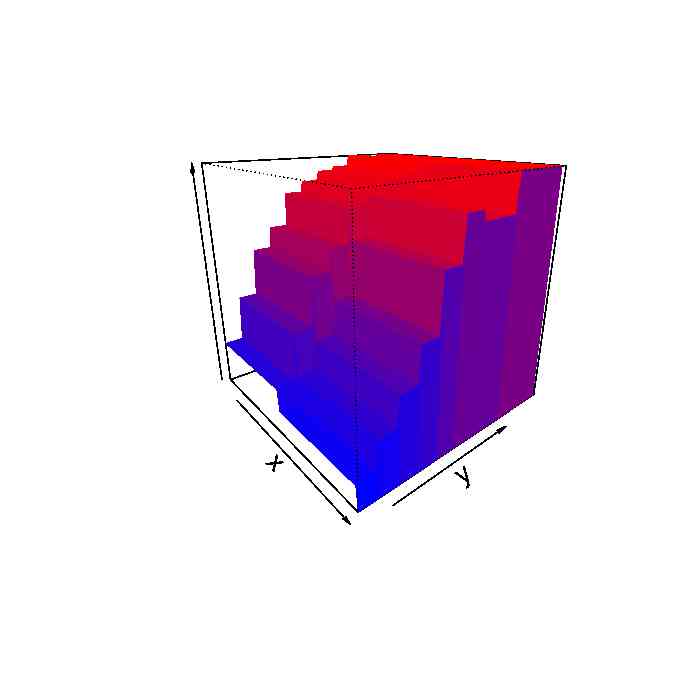} \\
\hspace{-1.5cm}\includegraphics[width=4cm,trim= 20mm 20mm 10mm 20mm]{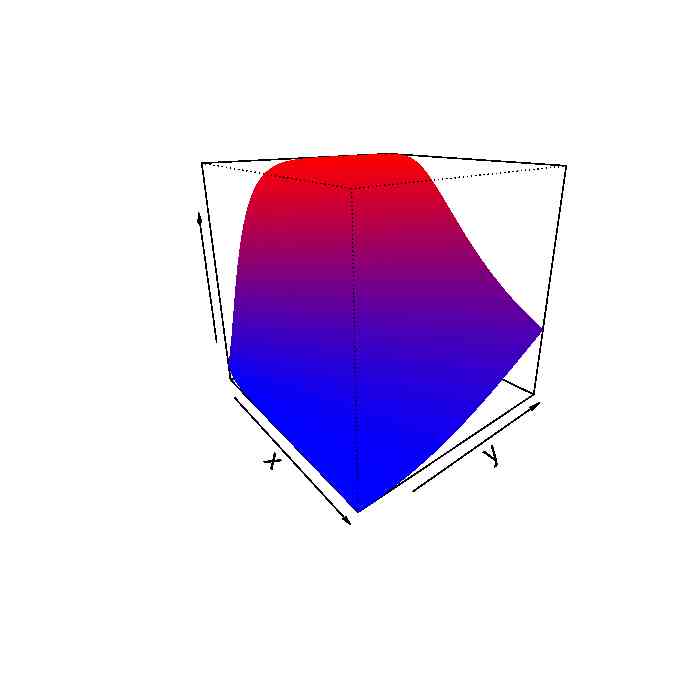} &
\includegraphics[width=4cm,trim= 20mm 20mm 10mm 20mm]{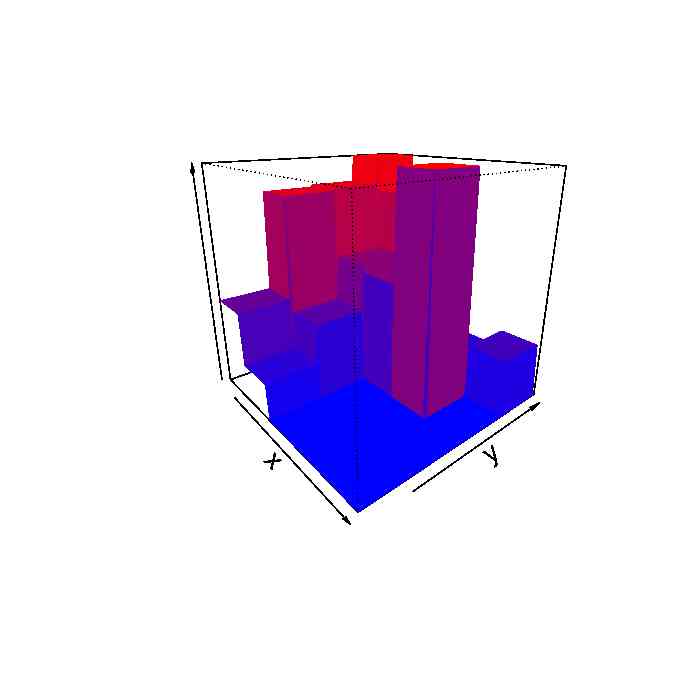} &
 \includegraphics[width=4cm,trim= 20mm 20mm 10mm 20mm]{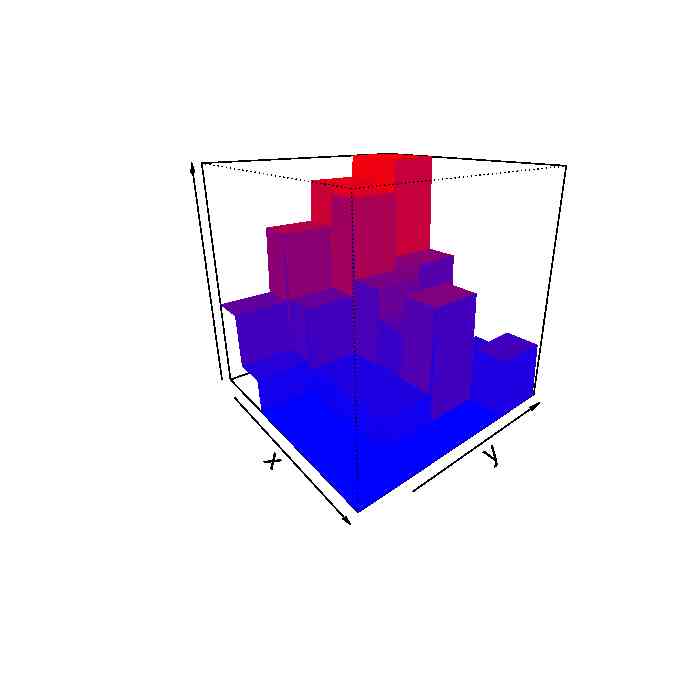} &
 \includegraphics[width=4cm,trim= 20mm 20mm 10mm 20mm]{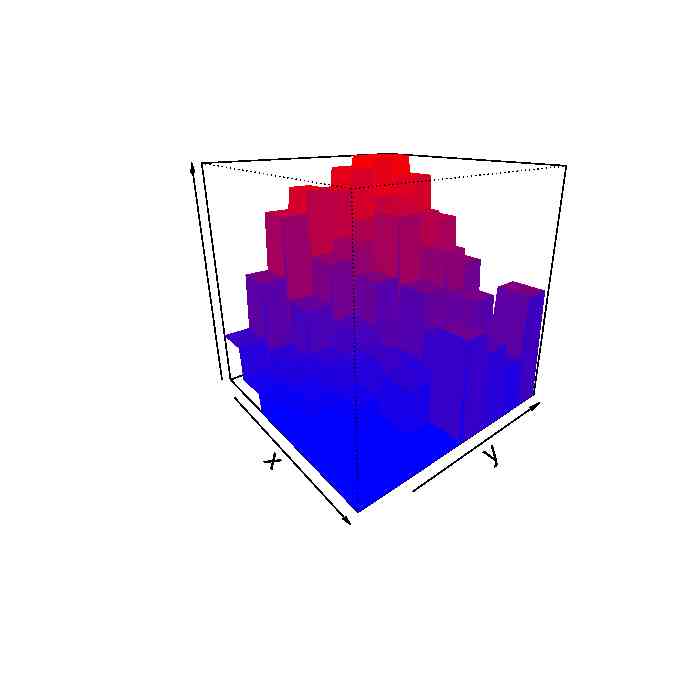}  \\
& $n=500$ & $n=1000$ & $ n=5000$  
\end{tabular}
\caption{\label{Fig-Fm} Conditional cumulative distribution function $F$ (first column) and model selection estimator $\tilde{F}_{\widehat{m}}$ for sample size $n=(500,1000,5000)$. First row: model 1; Second row: model 2.}
\end{figure}

\begin{figure}
\begin{tabular}{ccc}
\hspace{-1.2cm} \includegraphics[width=5cm,trim= 20mm 20mm 10mm 20mm]{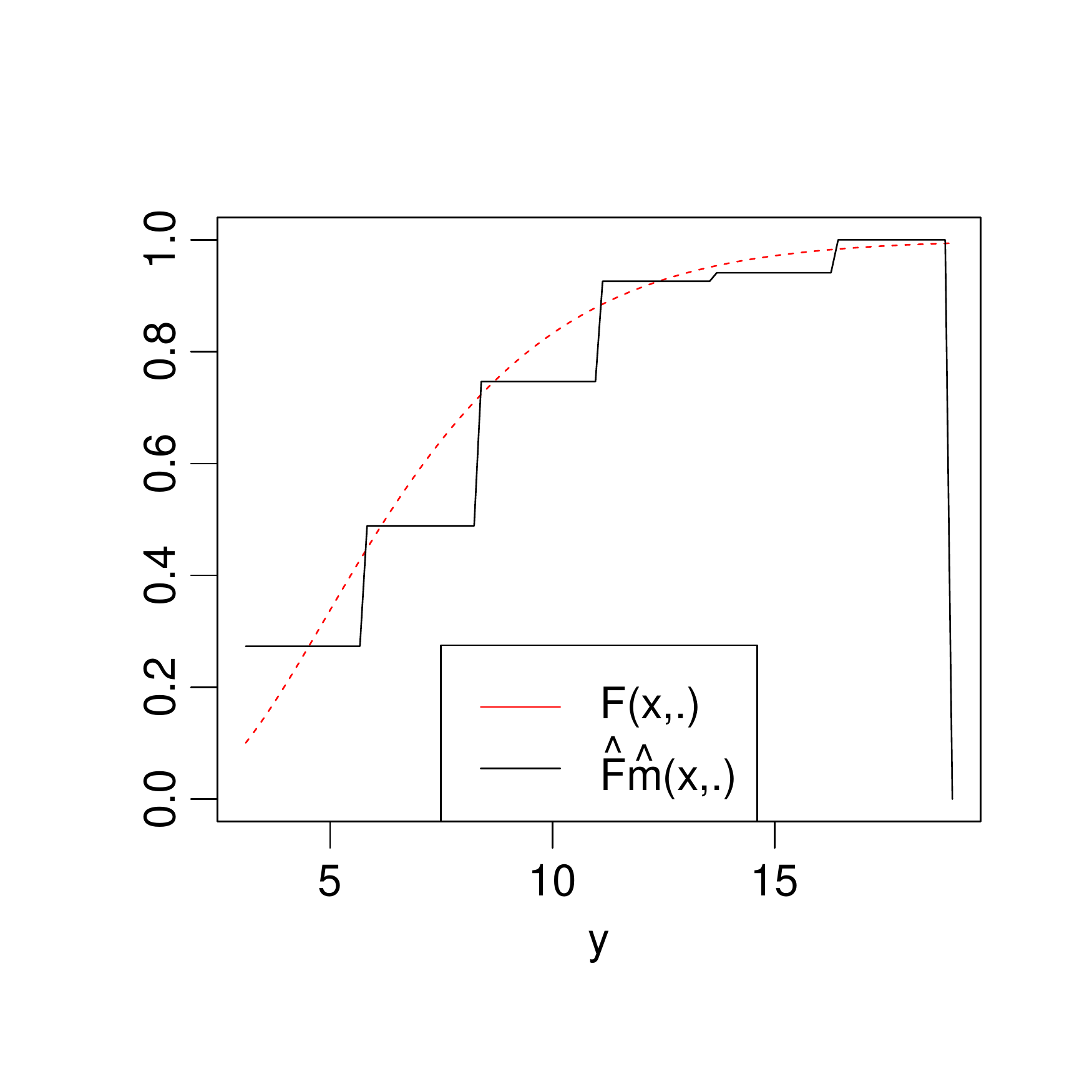} & 
\includegraphics[width=5cm,trim= 20mm 20mm 10mm 20mm]{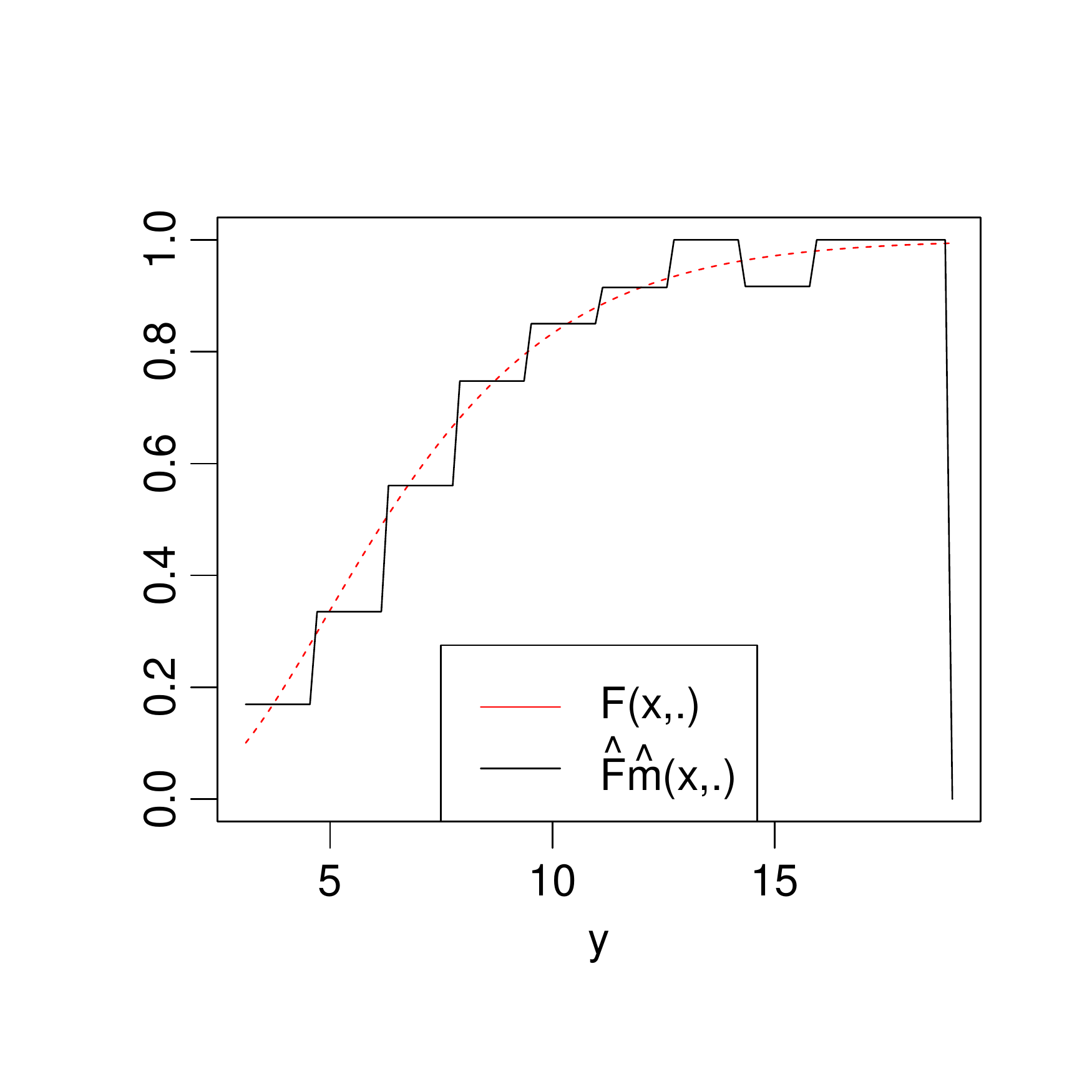} &
\includegraphics[width=5cm,trim= 20mm 20mm 10mm 20mm]{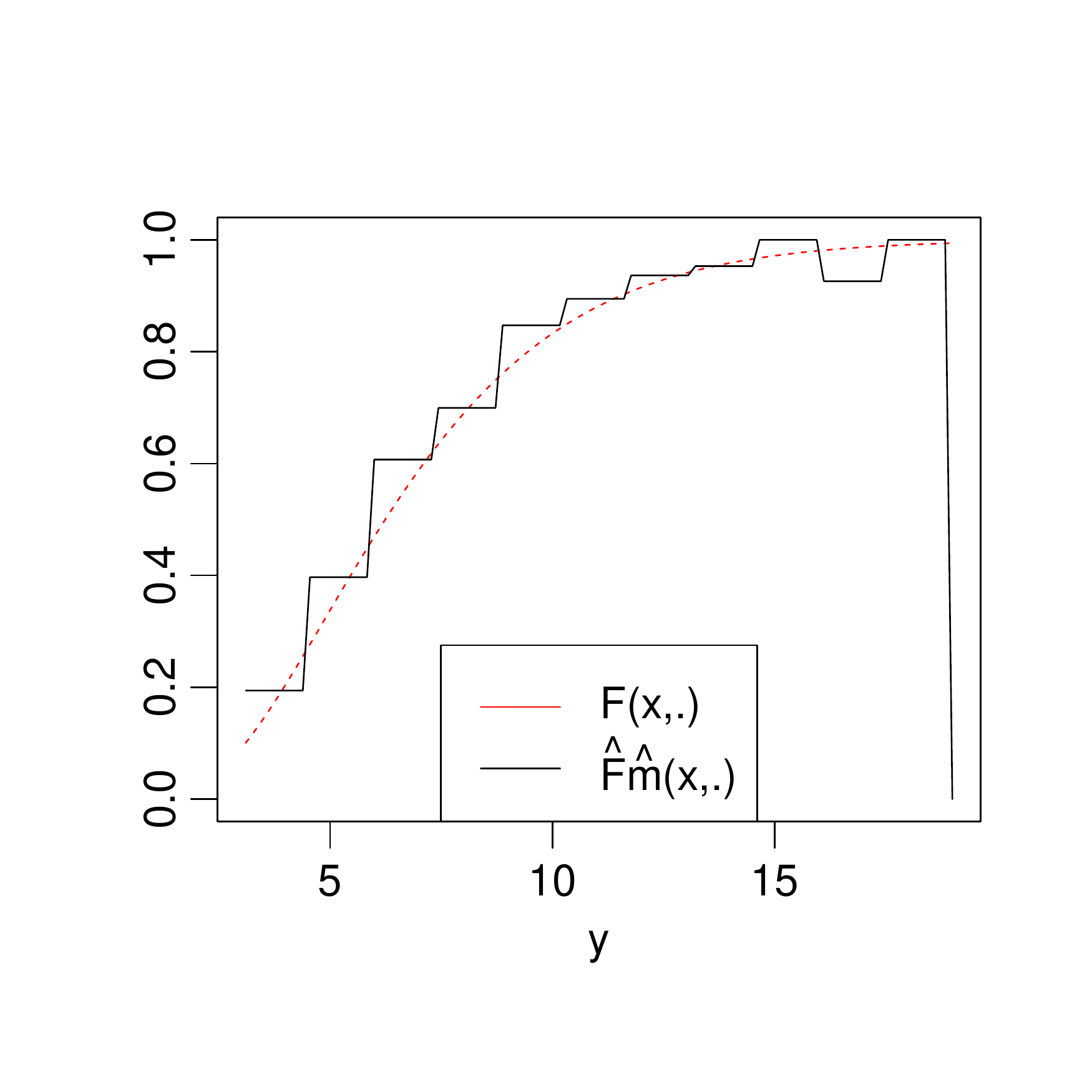} \\
\hspace{-1.2cm}\includegraphics[width=5cm,trim= 20mm 20mm 10mm 20mm]{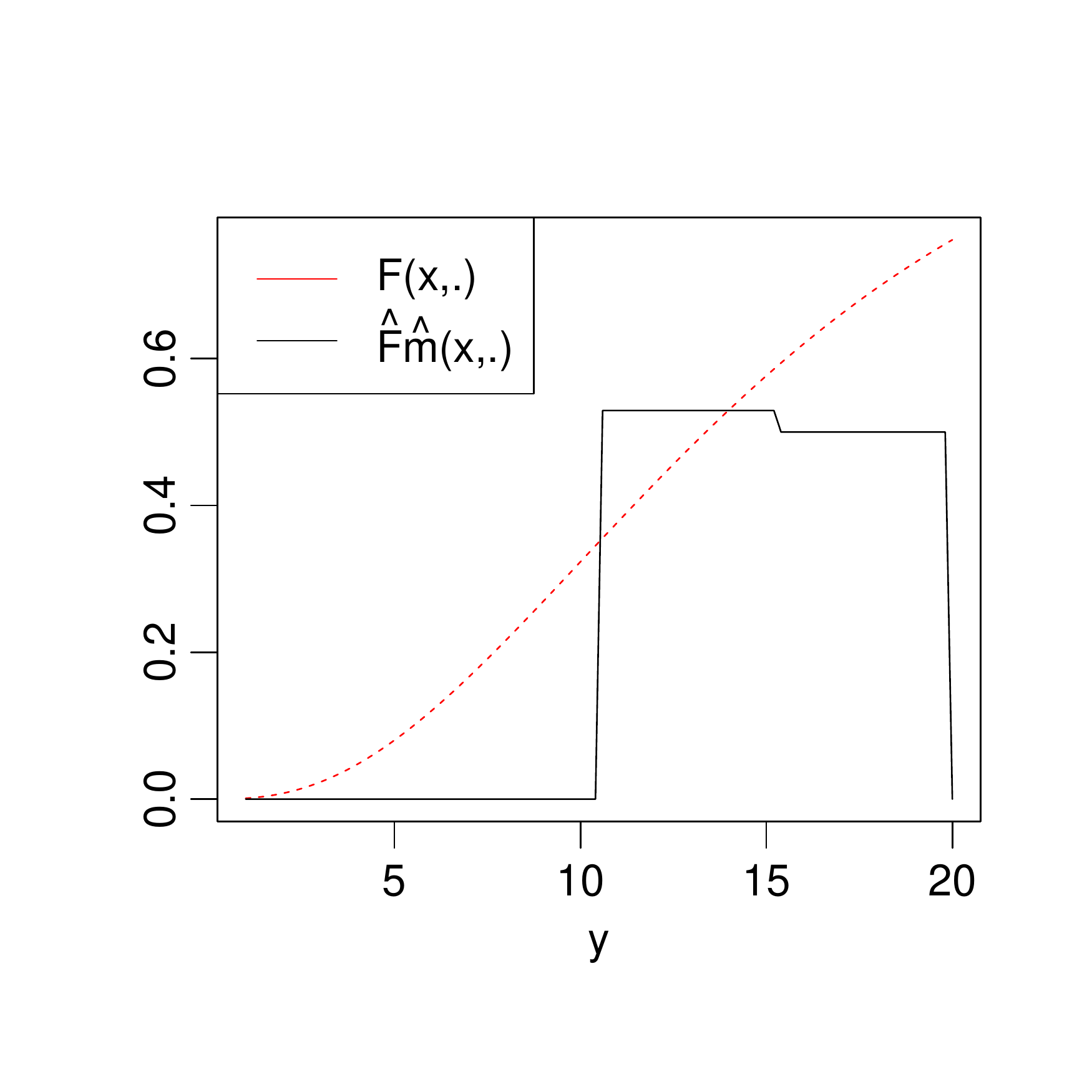} &
 \includegraphics[width=5cm,trim= 20mm 20mm 10mm 20mm]{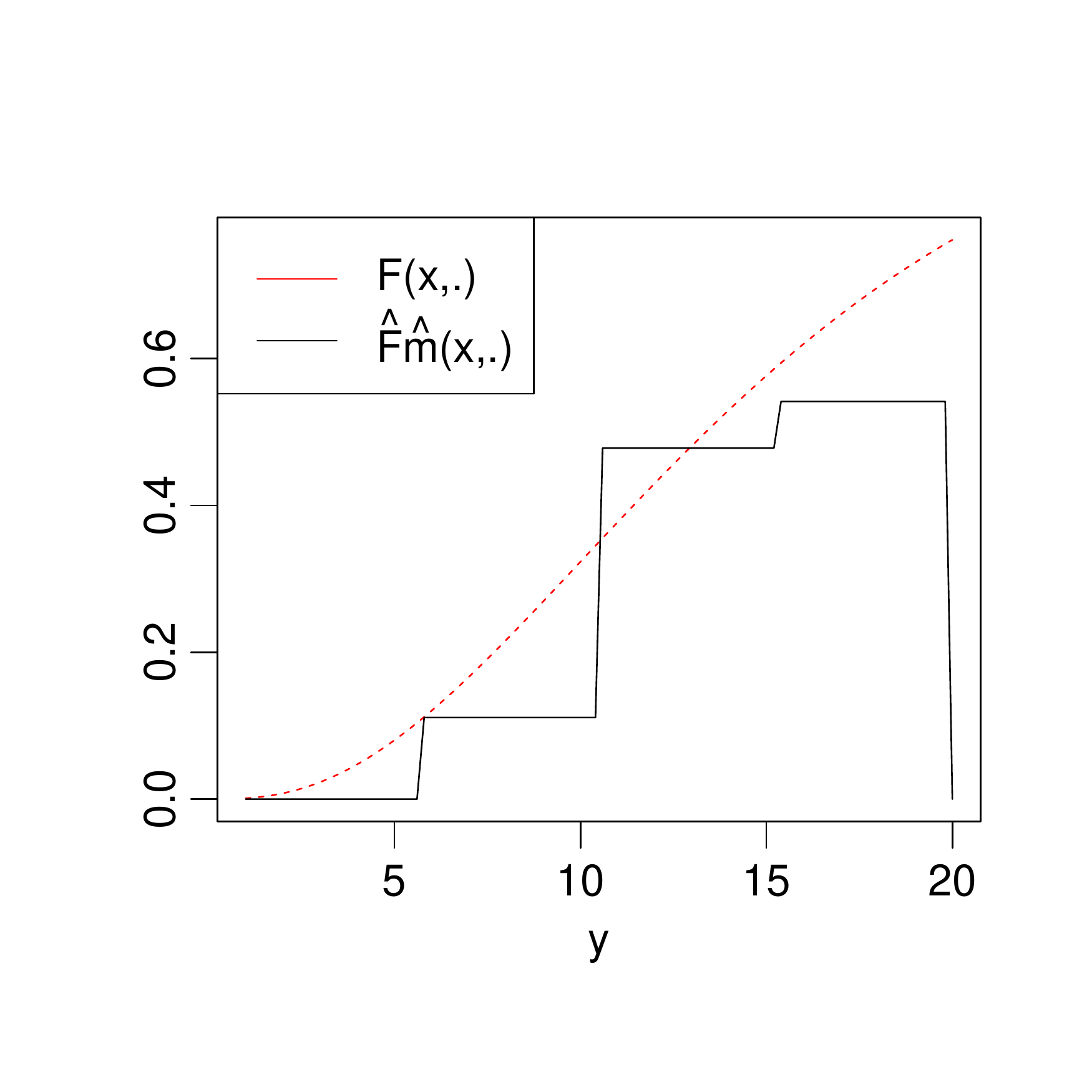}&
 \includegraphics[width=5cm,trim= 20mm 20mm 10mm 20mm]{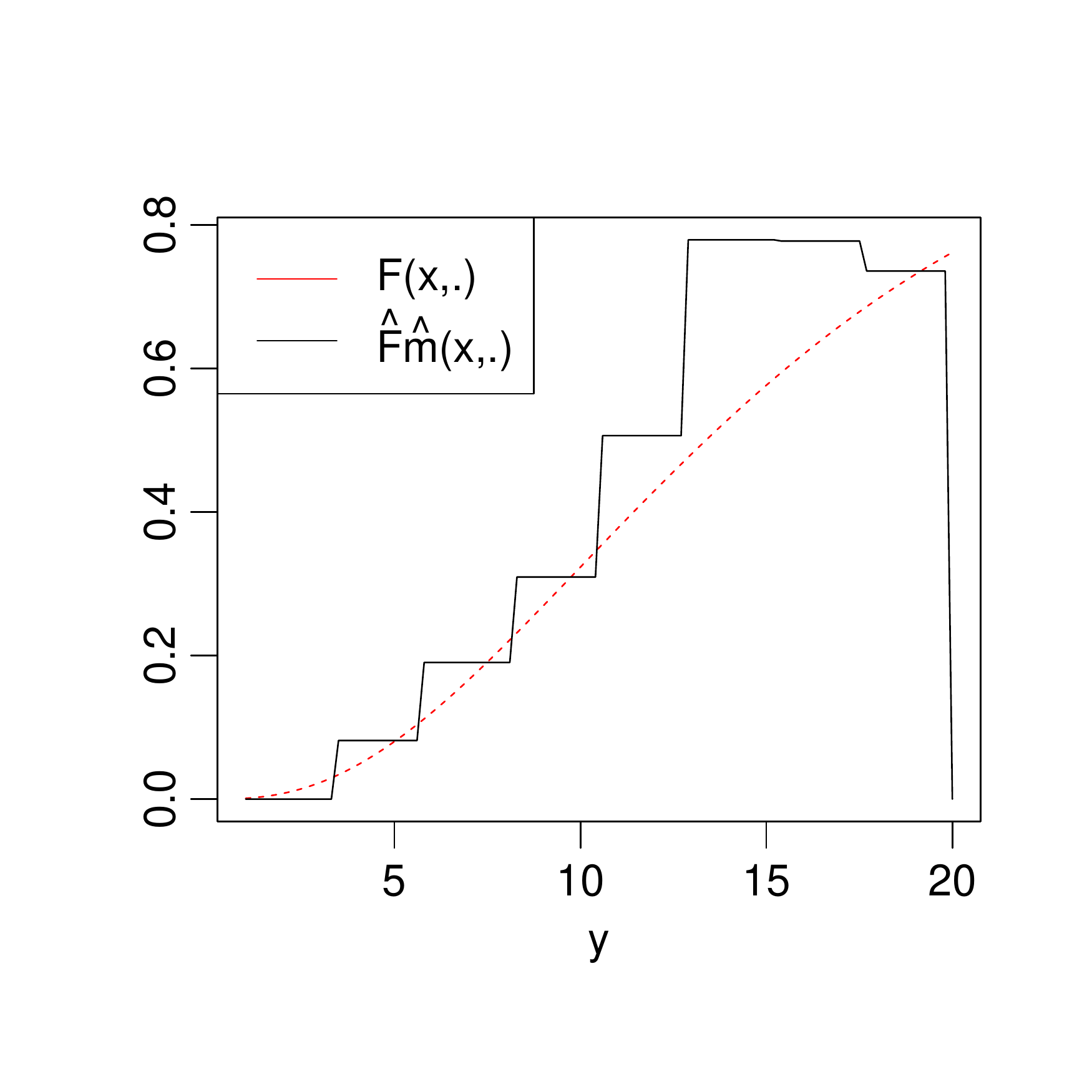} \\
 $n=500$ & $n=1000$ & $ n=5000$ 
\end{tabular}
\caption{\label{Fig-plotx} Model selection estimator $\tilde{F}_{\widehat{m}}$ (black solid line) and target function $F$ (red dotted line) for a fixed value of $x$. First row: model 1, $x=0.9$; second row: model 2, $x=5$.}
\end{figure}

\begin{figure}
\begin{tabular}{ccc}
\hspace{-1.2cm}\includegraphics[width=5cm,trim= 20mm 20mm 10mm 20mm]{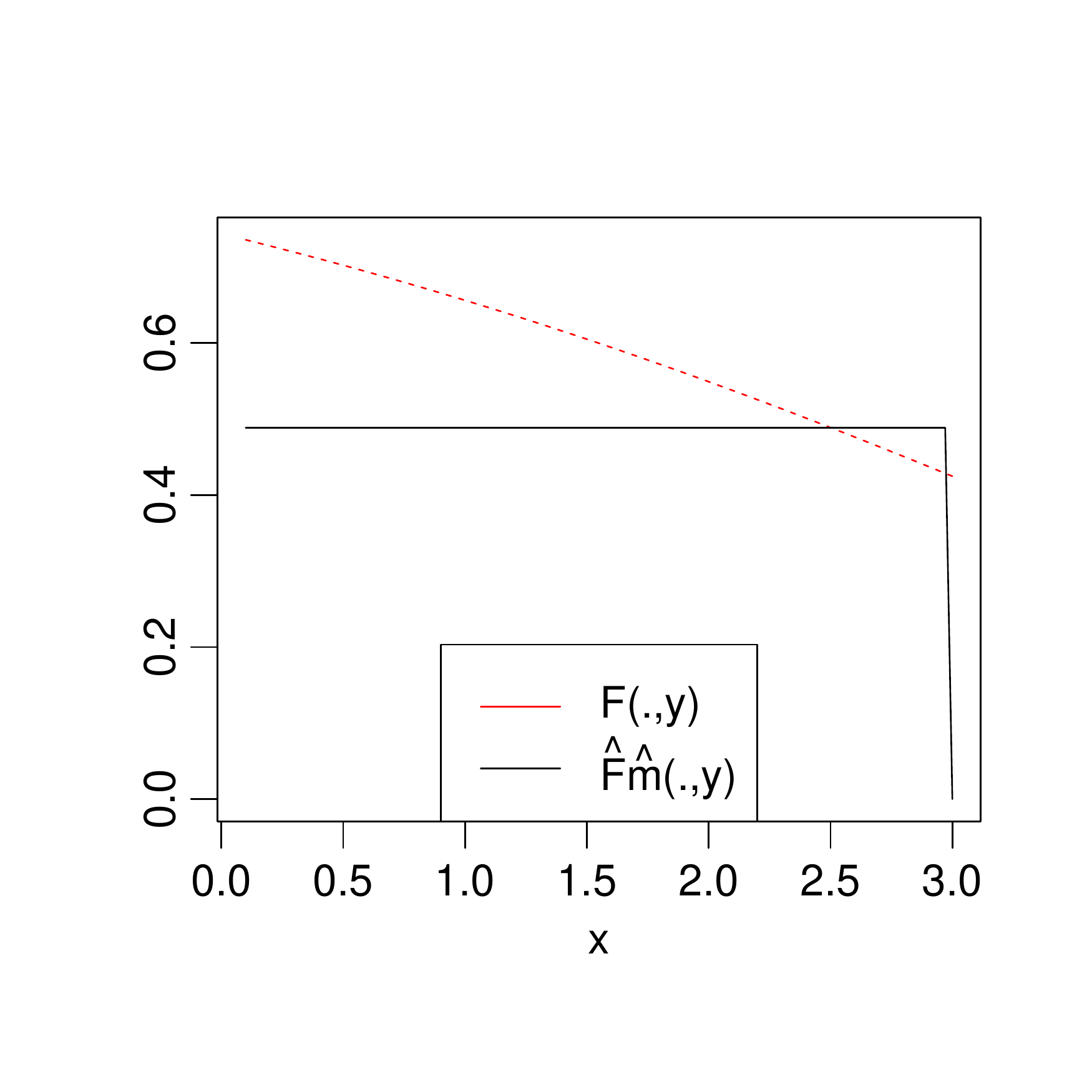} &
 \includegraphics[width=5cm,trim= 20mm 20mm 10mm 20mm]{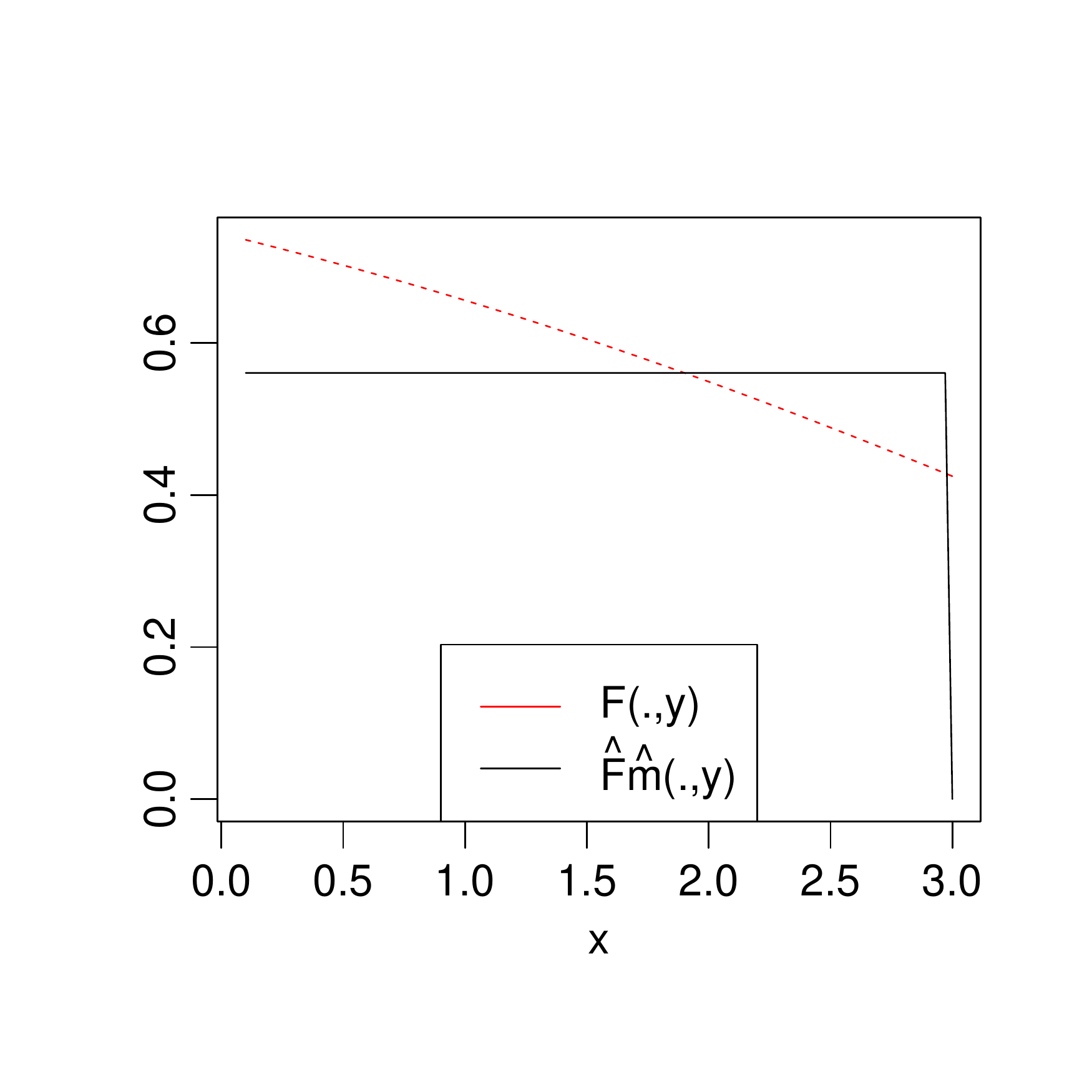} &
 \includegraphics[width=5cm,trim= 20mm 20mm 10mm 20mm]{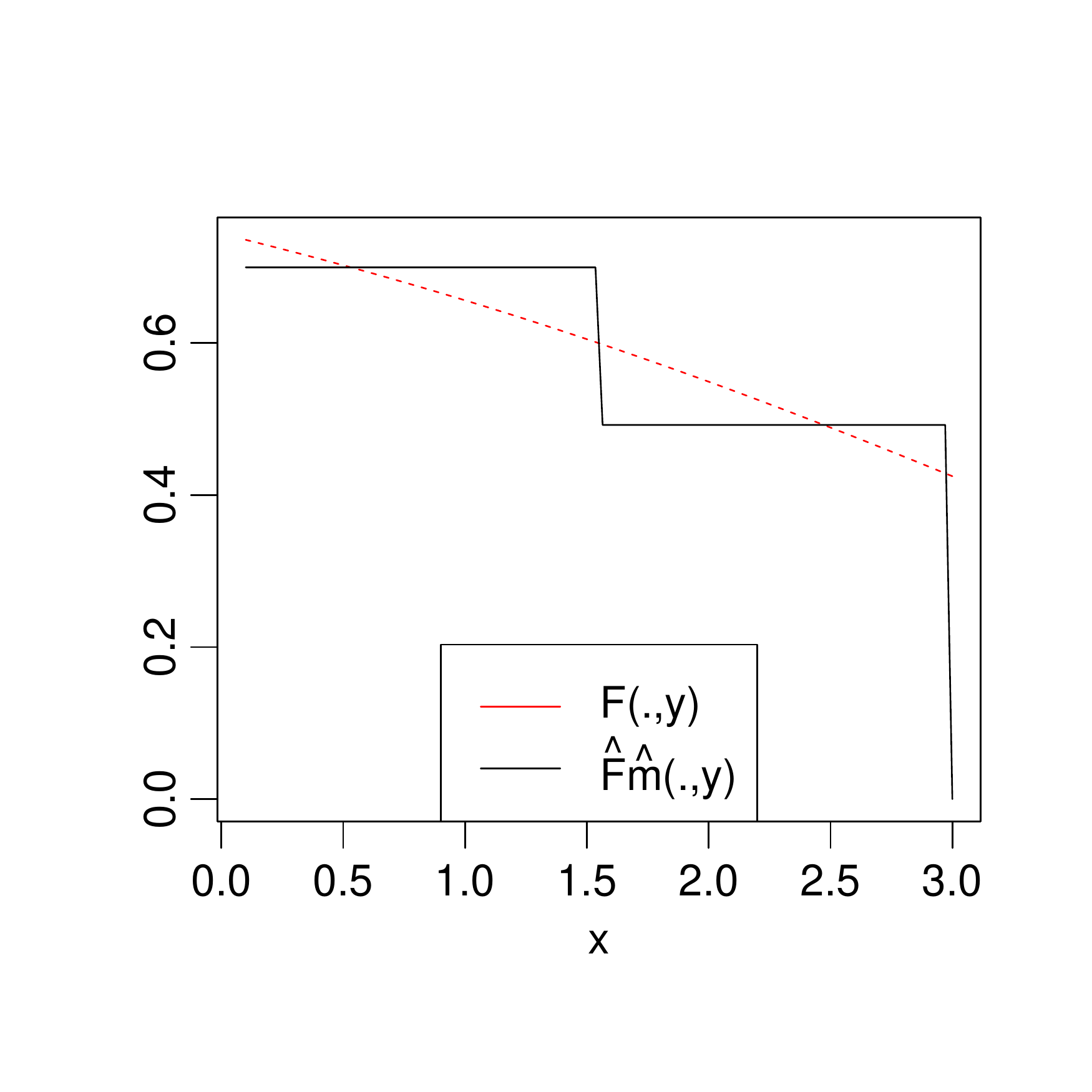} \\
\hspace{-1.2cm}\includegraphics[width=5cm,trim= 20mm 20mm 10mm 20mm]{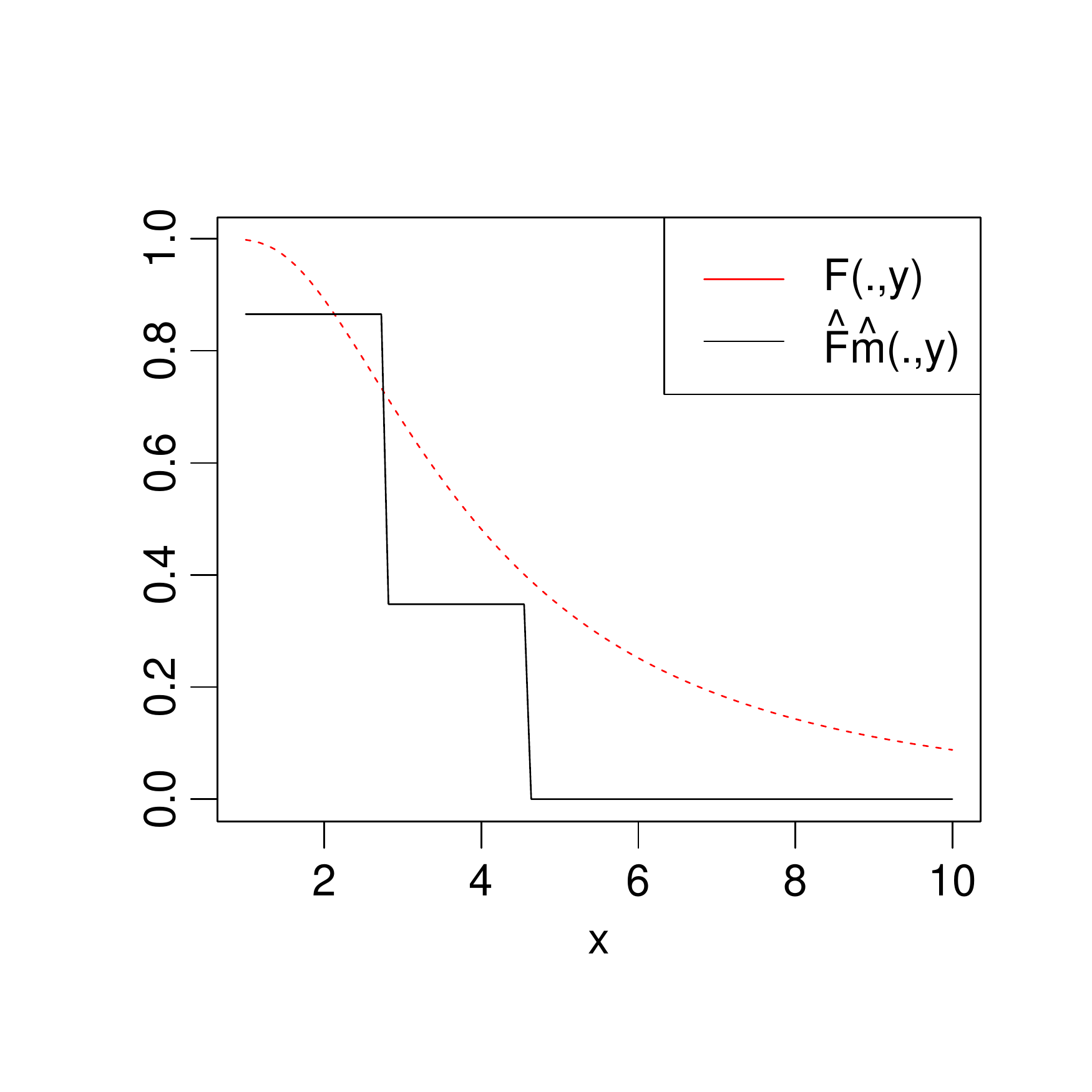} & 
\includegraphics[width=5cm,trim= 20mm 20mm 10mm 20mm]{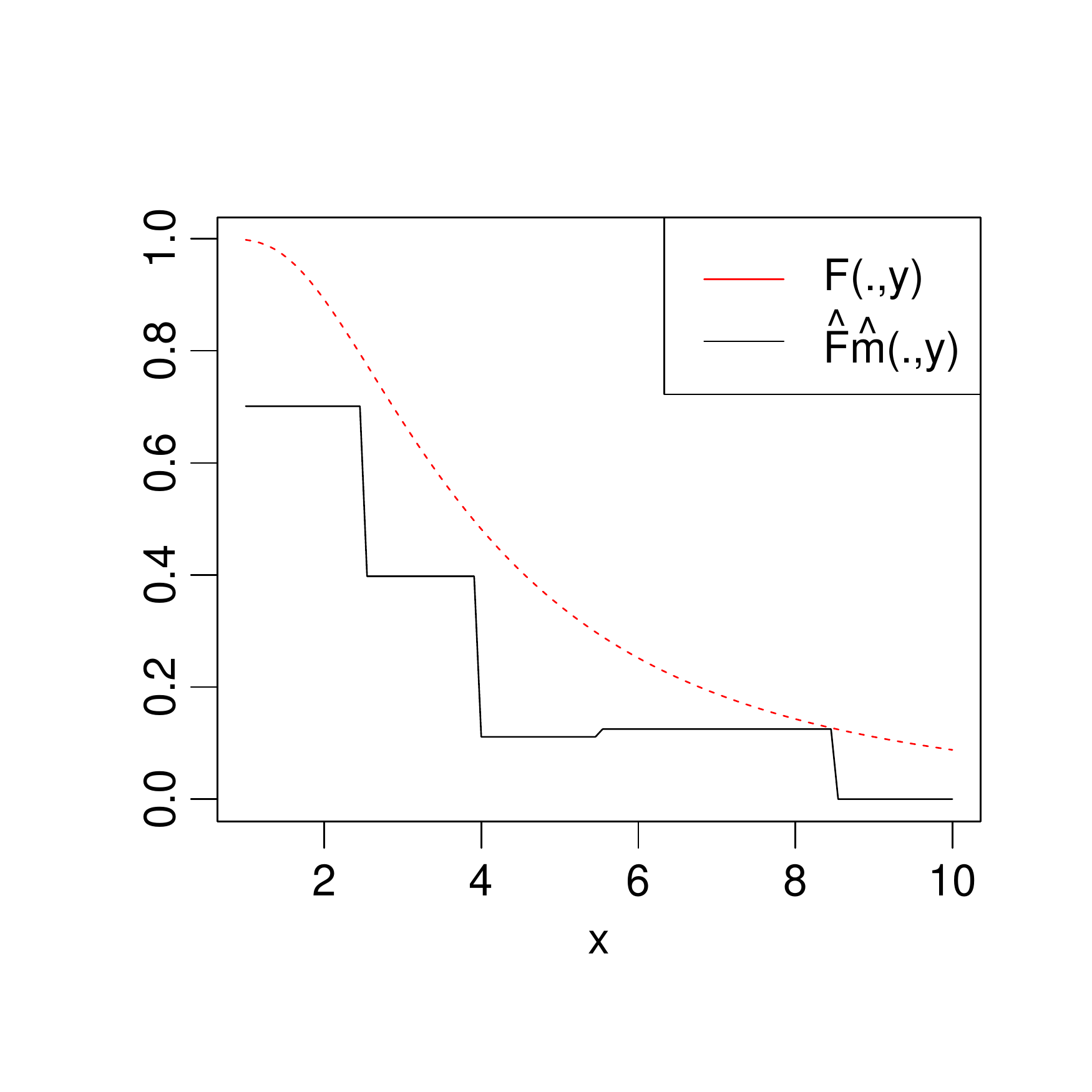} &
\includegraphics[width=5cm,trim= 20mm 20mm 10mm 20mm]{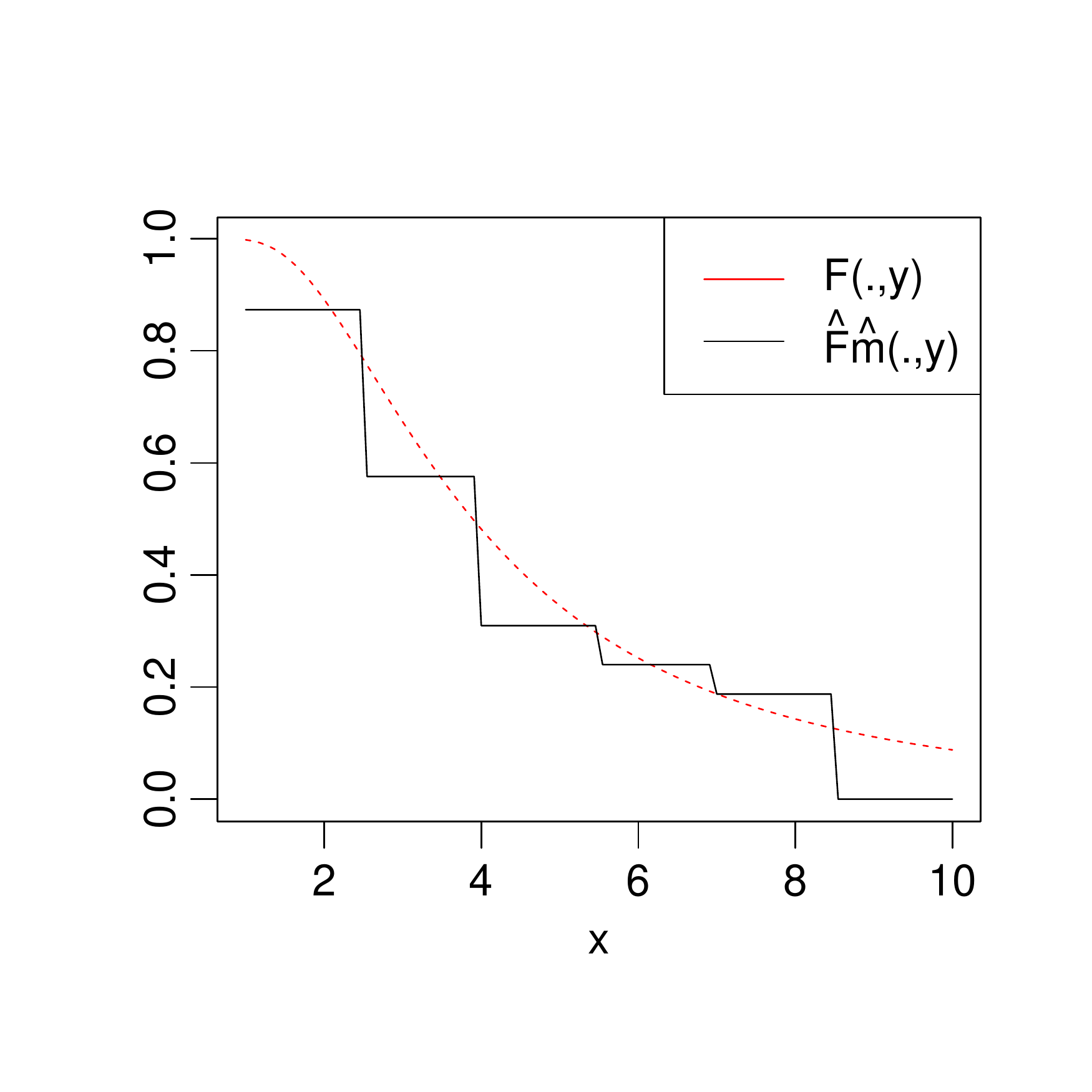} \\
 $n=500$ & $n=1000$ & $ n=5000$ 
\end{tabular}
\caption{\label{Fig-ploty} Model selection estimator $\tilde{F}_{\widehat{m}} $(black solid line) and target function $F$ (red dotted line) for a fixed value of $y$. First row: model , $y=7.5$; Second row: model 2, $y=10.5$.}
\end{figure}

For models 1 and 2, samples $(X_i,Y_i,T_i)_{i=1,\dots,n}$ of size $n=(500,1000,5000)$ are generated, and  the model selection estimators $\tilde{F}_{\widehat{m}}$ is computed from the sample  $(X_i,T_i, \delta_i=1\!\text{I} _{\{Y_i\leq T_i\}}) _{i=1,\dots,n} $. The results are displayed in Figure \ref{Fig-Fm}. Each row corresponds to a simulation model. The first column presents the target function $F$ and columns 2 to 4 present the estimators for sample size $n=(500,1000,5000)$. The same functions are plotted in Figures \ref{Fig-plotx} and \ref{Fig-ploty} for a fixed $x$ and $y$ respectively. As expected, the estimation gets better when the sample size increases, but we notice that a large sample size is required in order to get correct estimation. Indeed, on the one hand,  the nature of the current status data provides less precise information about the variable of interest than uncensored of right-censored data, and on the other hand, the transition from uni- to bi-dimensional framework induces a loss in estimation quality for a given sample size.

\subsection{Impact of the distance between the densities of $Y$ and $T$}

In a right-censoring framework, the rate of censoring, defined as the expected proportion of observations that are censored, is a well-known factor which impacts the quality of estimation:  higher rate of censoring leads to poorer estimation. Indeed as the rate of censoring decreases, the proportion of survival times actually observed increases and the estimation of the  survival time distribution gets better.
In the current status data framework, the rate of censoring can not be defined since the survival time is never directly observed. In alternative, the distance between the densities of $Y$ and $T$ appears as a relevant parameter that impacts the quality of estimation. 
Let us heuristically explain this idea. 

First of all, the estimation of $F$ is more accurate if the observation times are concentrated in the areas where $F$ shows the highest variations. Besides, for a given value of the variable $X$, the highest variations of the c.d.f. $F$ correspond to the largest values of the conditional density $f_{Y|X}$. Therefore, the estimation of $F$ improves as the densities of the observation time and survival time get closer, which may be quantified by the distance between the marginal densities $f_Y$ and $f_T$.
As a particular case, consider the situation when the supports of $f_Y$ and $f_T$ are disconnected: $\delta _i$ has the same value, either 0 or 1, for all $i$; thus the observation set $(X_i,T_i,\delta _i)_{i=1\dots,n} $ do not provide information about the distribution of $Y$.  We enhance this heuristic on numerical examples. \\

\begin{figure}
\begin{tabular}{cccc}
\hspace{-1.5cm}\includegraphics[width=4cm,trim= 10mm 20mm 10mm 20mm]{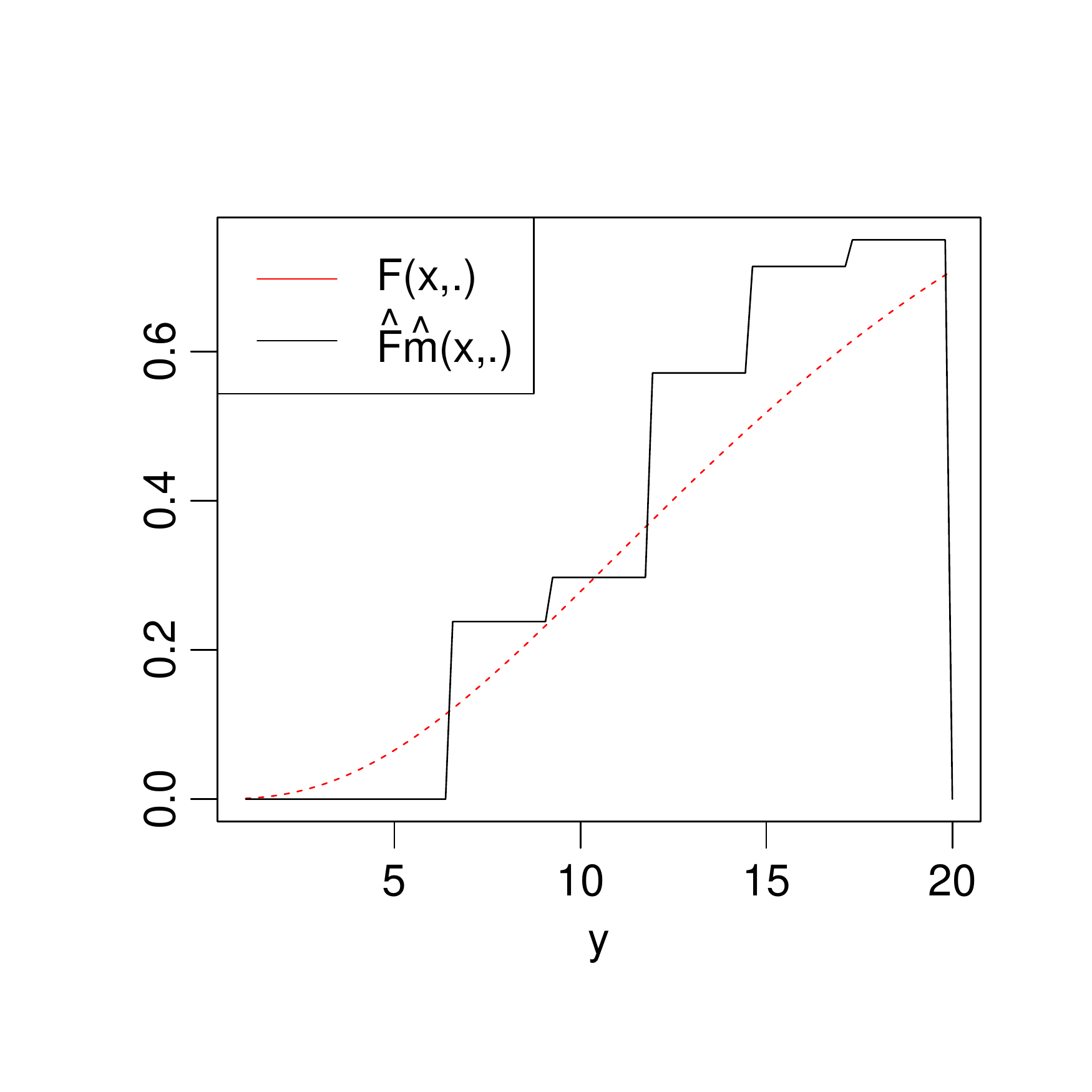} &
 \includegraphics[width=4cm,trim= 10mm 20mm 10mm 20mm]{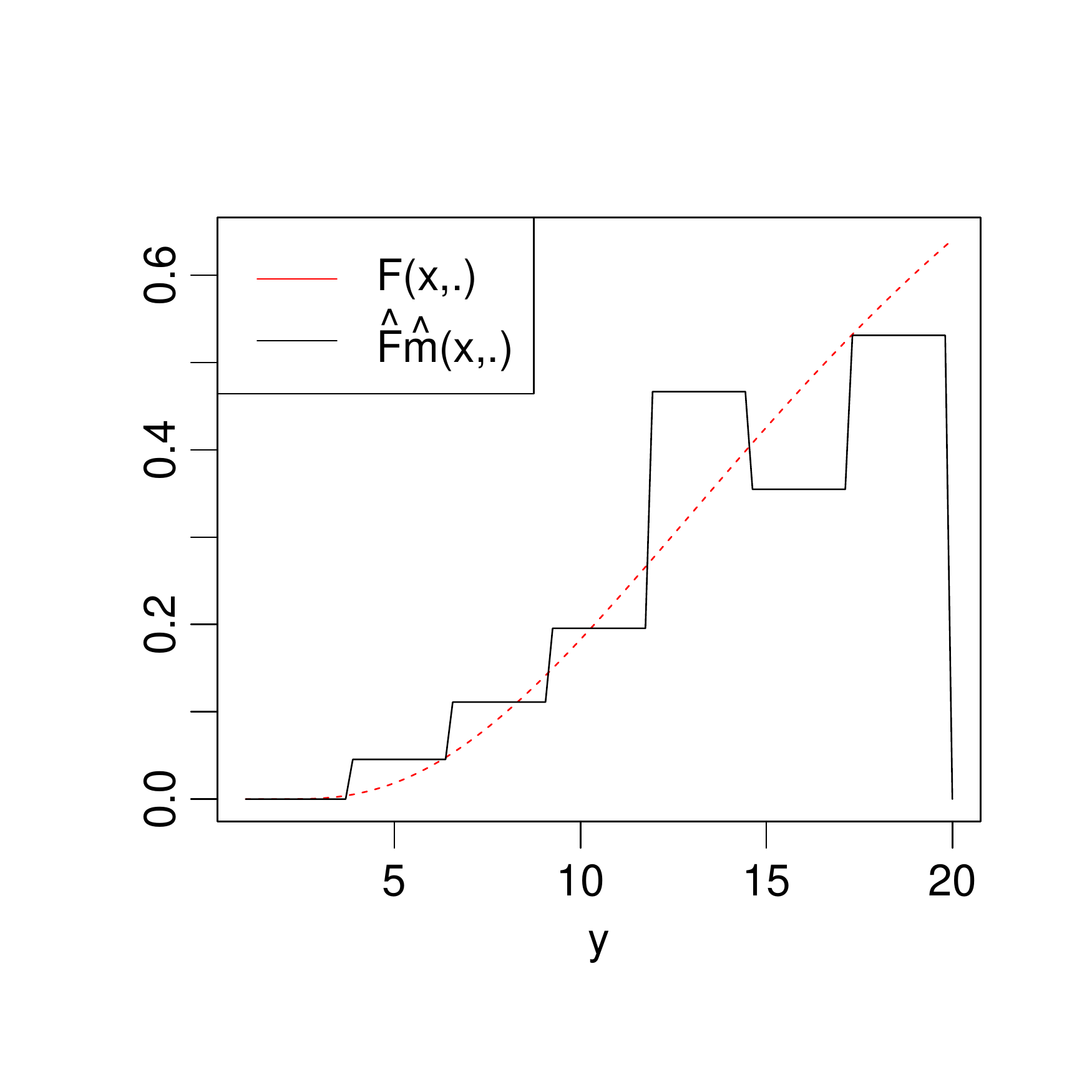} &
 \includegraphics[width=4cm,trim= 10mm 20mm 10mm 20mm]{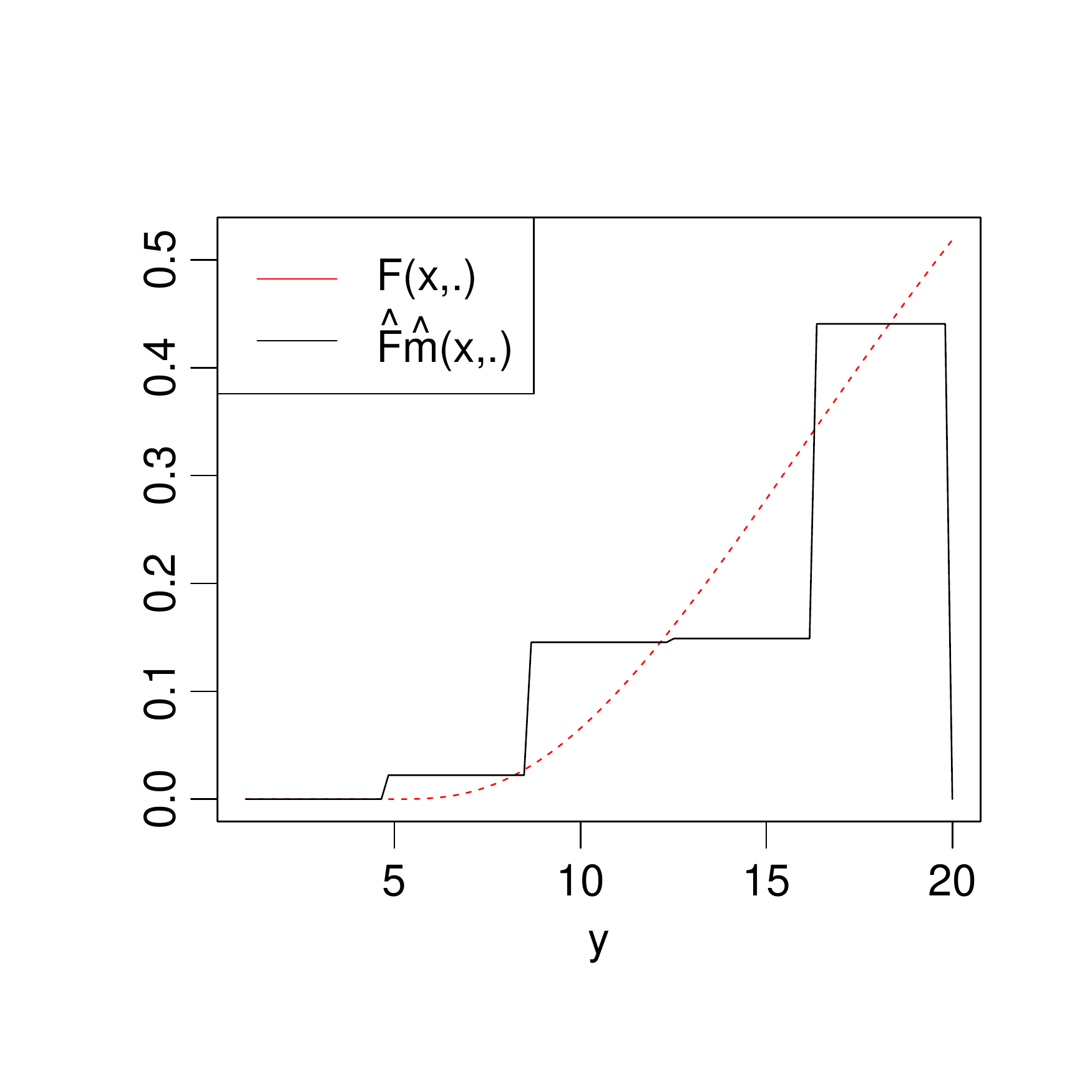}  &
  \includegraphics[width=4cm,trim= 10mm 20mm 10mm 20mm]{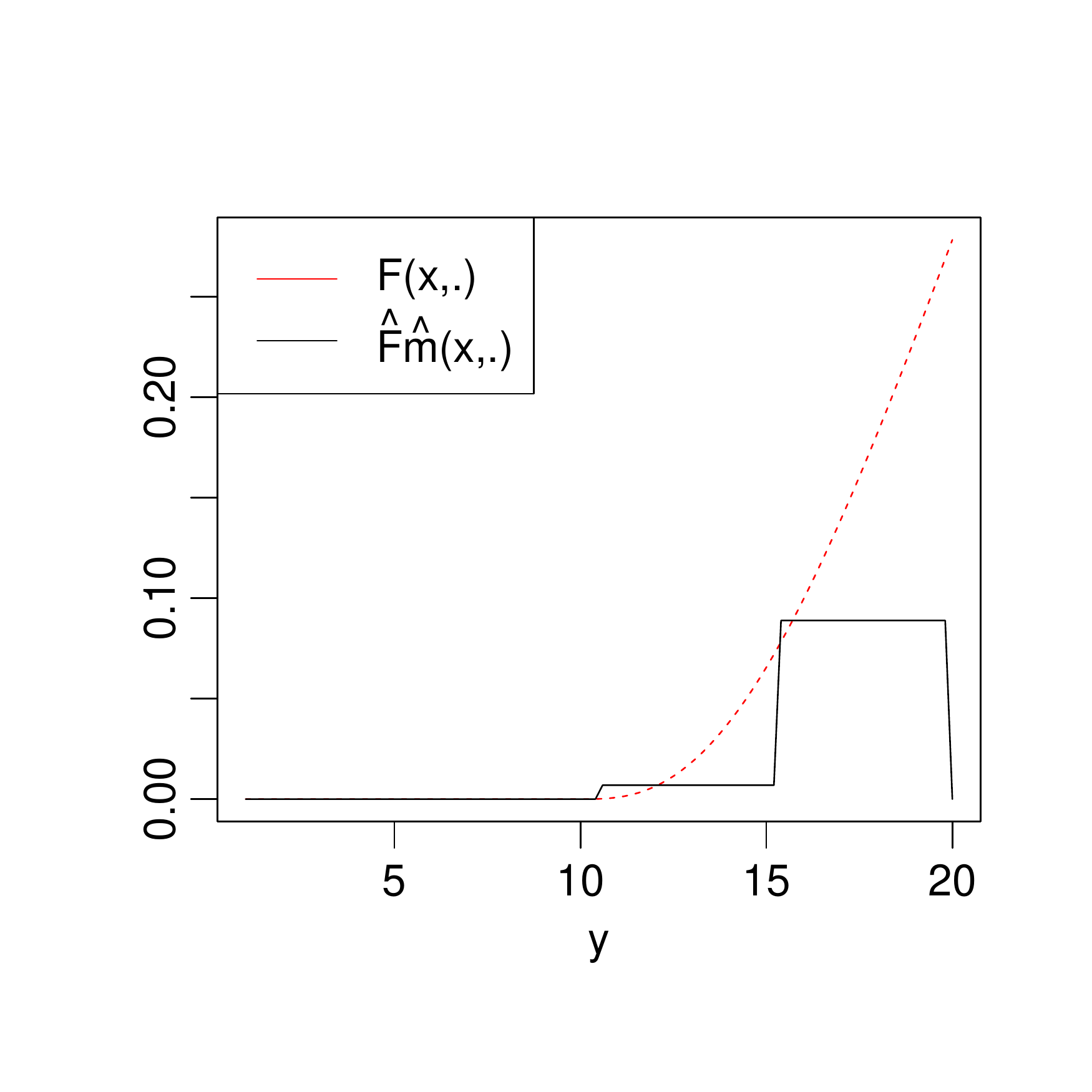}\\
$a=0$& $a=2$ & $a=5$ & $a=10$ 
\end{tabular}
\caption{\label{Fig-off} Results on simulation from Model 2b with offset $a =(0,2,5,10)$; model selection estimator (black solid line) and target function $F$ (red dotted line) for a fixed value of $x$, $x=5.5$. .}
\end{figure}

\noindent \textbf{Model 2b. } We consider the same distribution as Model 2 except that an offset $a$ is added to $Y$:

$$\begin{array}{l}
\\
\text{(Mod 2b)}\\
\\
\end{array}\left\{\begin{array}{l}
X\sim \Gamma (k=1.5,\theta=2) \\
Y = a+ X \times \varepsilon \quad \text{with} \quad  \varepsilon \sim \Gamma (k=3,\theta=1) \\
T = X\times \varepsilon' \quad \text{with} \quad  \varepsilon ' \sim \Gamma (k=3,\theta=1) \\
\end{array}
\right.$$
 The conditional cumulative distribution function $F(x,y) = \phi(y/x, k=3,\theta=1)$ is estimated on the set $A=[1,10]\times [1,20]$. The density $f_{(X,T)}$ is lower bounded by $h_0=10^{-7}$ on A. We consider the $L1$-distance between the densities of $Y$ and $T$:

\begin{eqnarray*}
 \text{dist}(a) &=& \| f_Y -f_T \| _{L^1} = \int _{x=0}^{+\infty} \int_{u=0}^{+\infty} \left| f_{Y|X}(x,u) -f_{T|X}(x,u)  \right| f_X(x) dx du \\
&=&
 \int _{x=0}^{+\infty} \int_{u=a}^{+\infty} \left| \phi \left(\frac{u-a}{x},3,1 \right) - \phi \left(\frac{u}{x},3,1 \right) \right| \phi (x,1.5,2) dx du 
\end{eqnarray*}
The values of $\text{dist} (a)$ for $a = (0,2,5,10)$ are displayed in Table \ref{tab-dista}.
For each value of $a$, a sample $(X_i,Y_i,T_i)_{i=1,\dots,n}$ of size $n=3000$ is generated from Model 2b, and the model selection estimator is computed on the set $A=[1,10]\times [1,20]$. The results are presented in Figure \ref{Fig-off}. A degradation in the estimation is observed as the distance between the densities of $Y$ and $T$ increases. In particular, for $a=10$, dist($a$) gets close to 2, which indicates that the supports of the distribution of $Y$ and $T$ hardly overlap, and the estimation is very poor despite a large sample size. 

\begin{table}\label{tab-dista}
\begin{center}
\begin{tabular}{ccccc} \hline
\rule{0pt}{15pt} \vspace{0.2cm} $a$ & 0  & 2 & 5 & 10 \\ \hline
 \rule{0pt}{15pt} \vspace{0.2cm} dist($a$) & 0 & 0.63 & 1.12 & 1.54  \\\hline
\end{tabular}
\vspace{0.2cm}
\end{center}
\caption{ \label{tab-dista} $L^1$-distance between $f_Y$ and $f_T$ for an offset $a$.} 
\end{table}

\subsection{Increasing rearrangement}
\begin{figure}
\begin{tabular}{ccc}
\hspace{-1.2cm}\includegraphics[width=5cm,trim= 20mm 20mm 10mm 20mm]{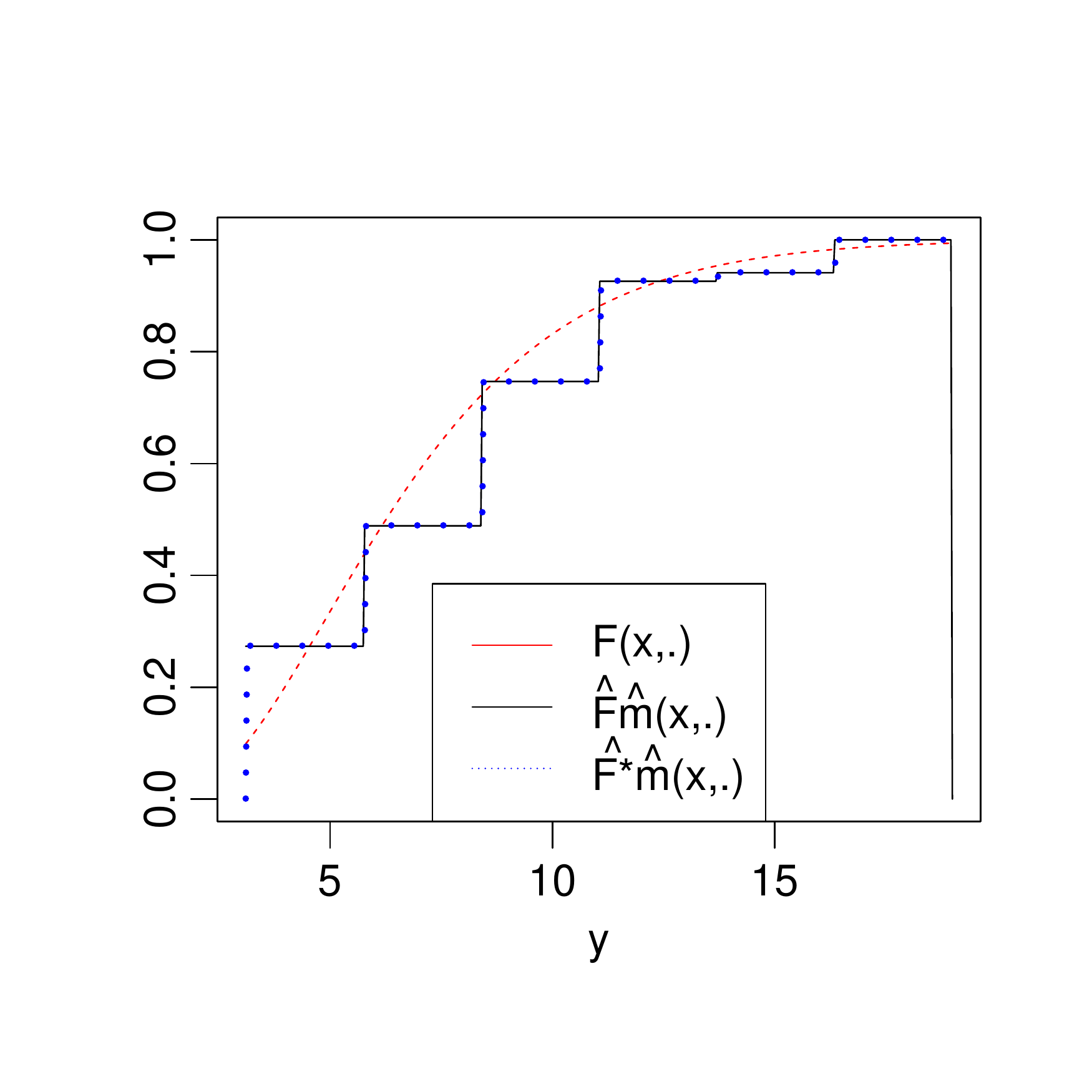} & 
\includegraphics[width=5cm,trim= 20mm 20mm 10mm 20mm]{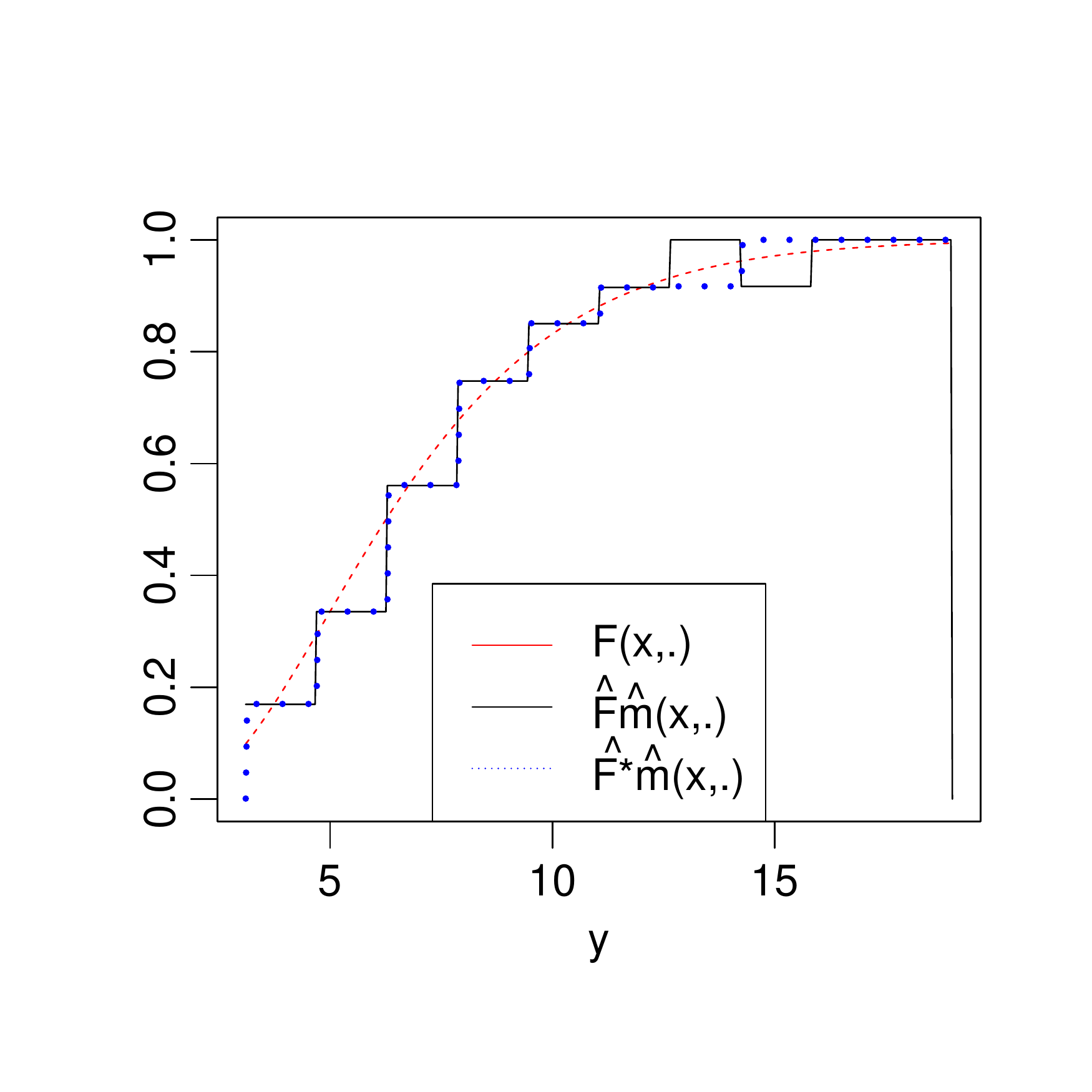} &
\includegraphics[width=5cm,trim= 20mm 20mm 10mm 20mm]{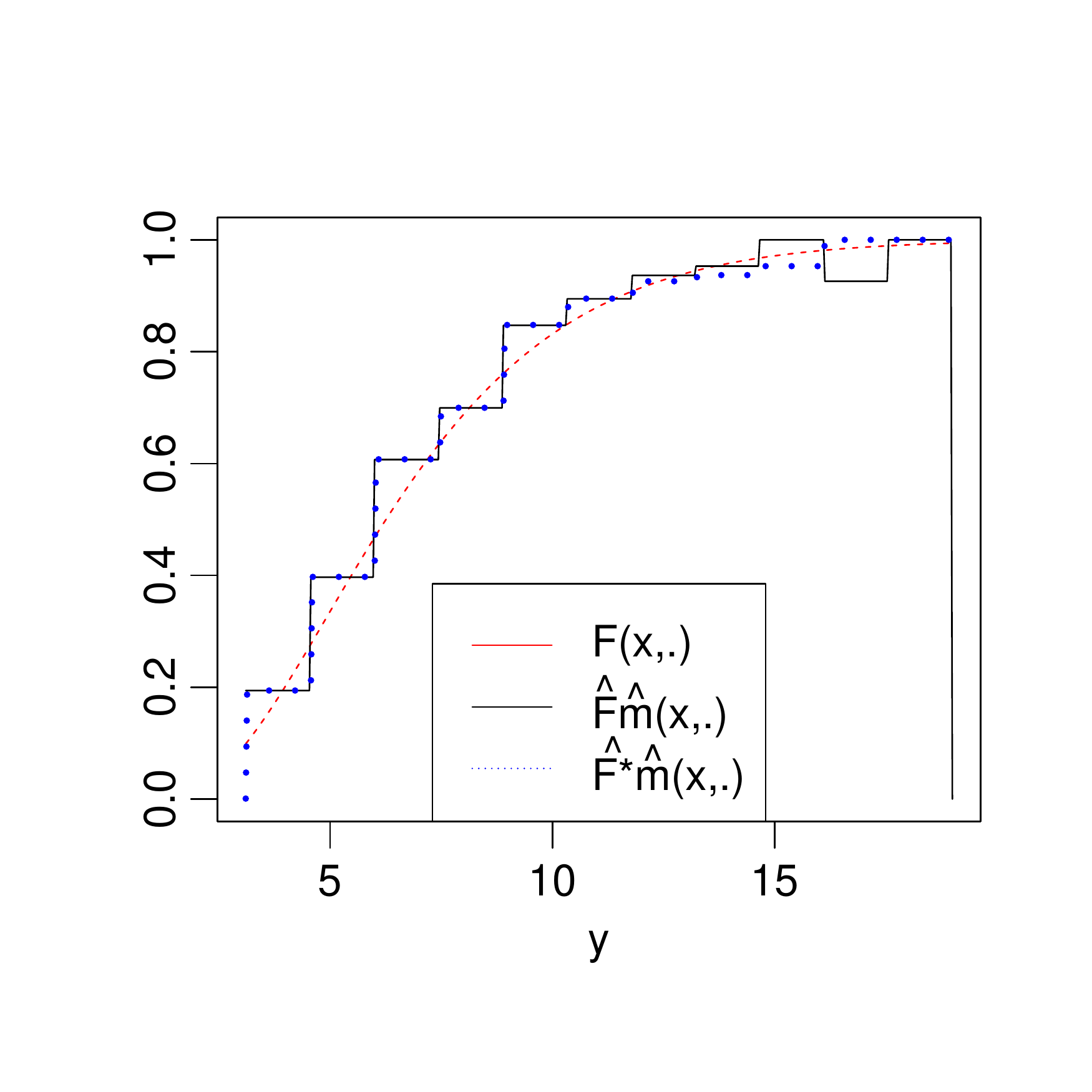} \\
\hspace{-1.2cm}\includegraphics[width=5cm,trim= 20mm 20mm 10mm 20mm]{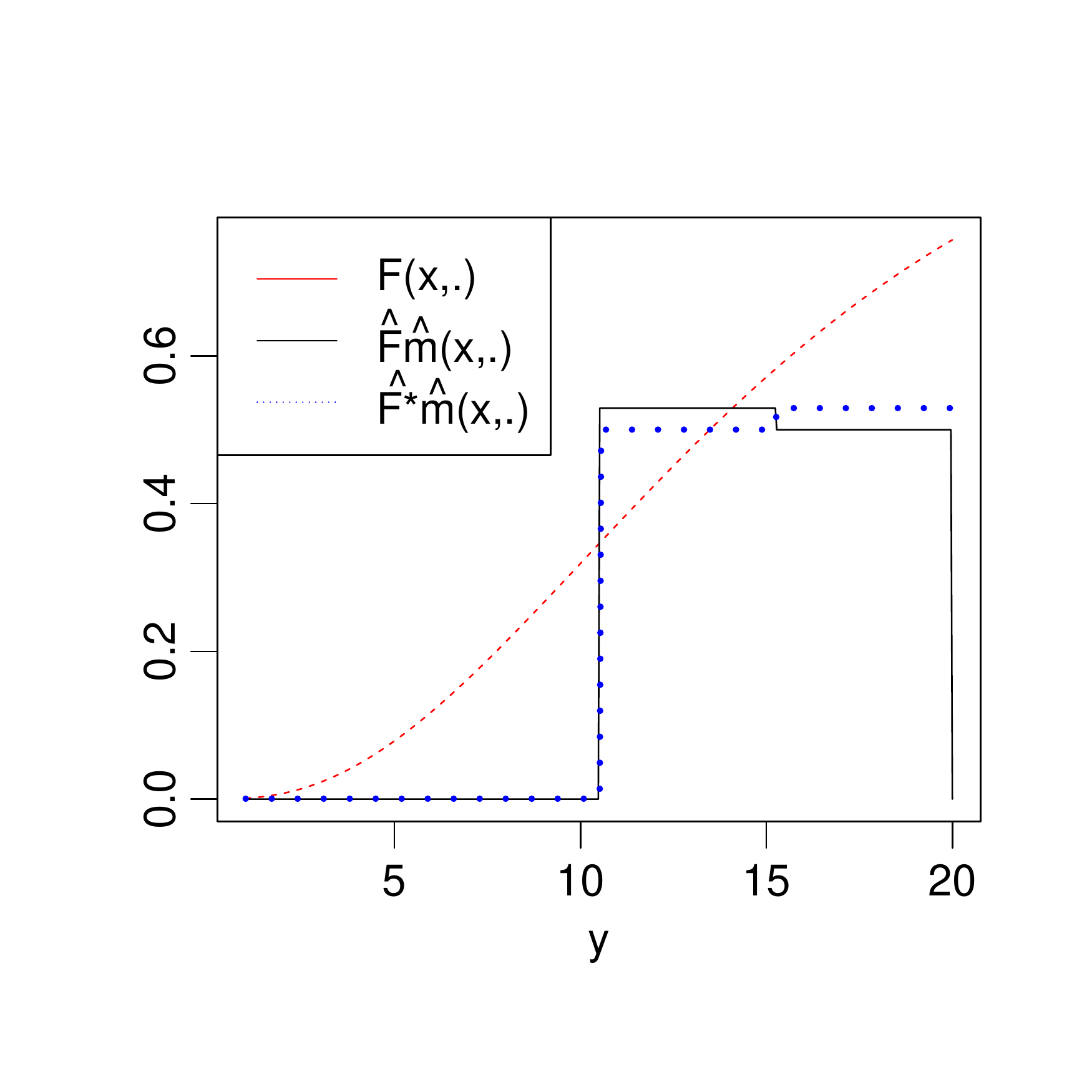} &
 \includegraphics[width=5cm,trim= 20mm 20mm 10mm 20mm]{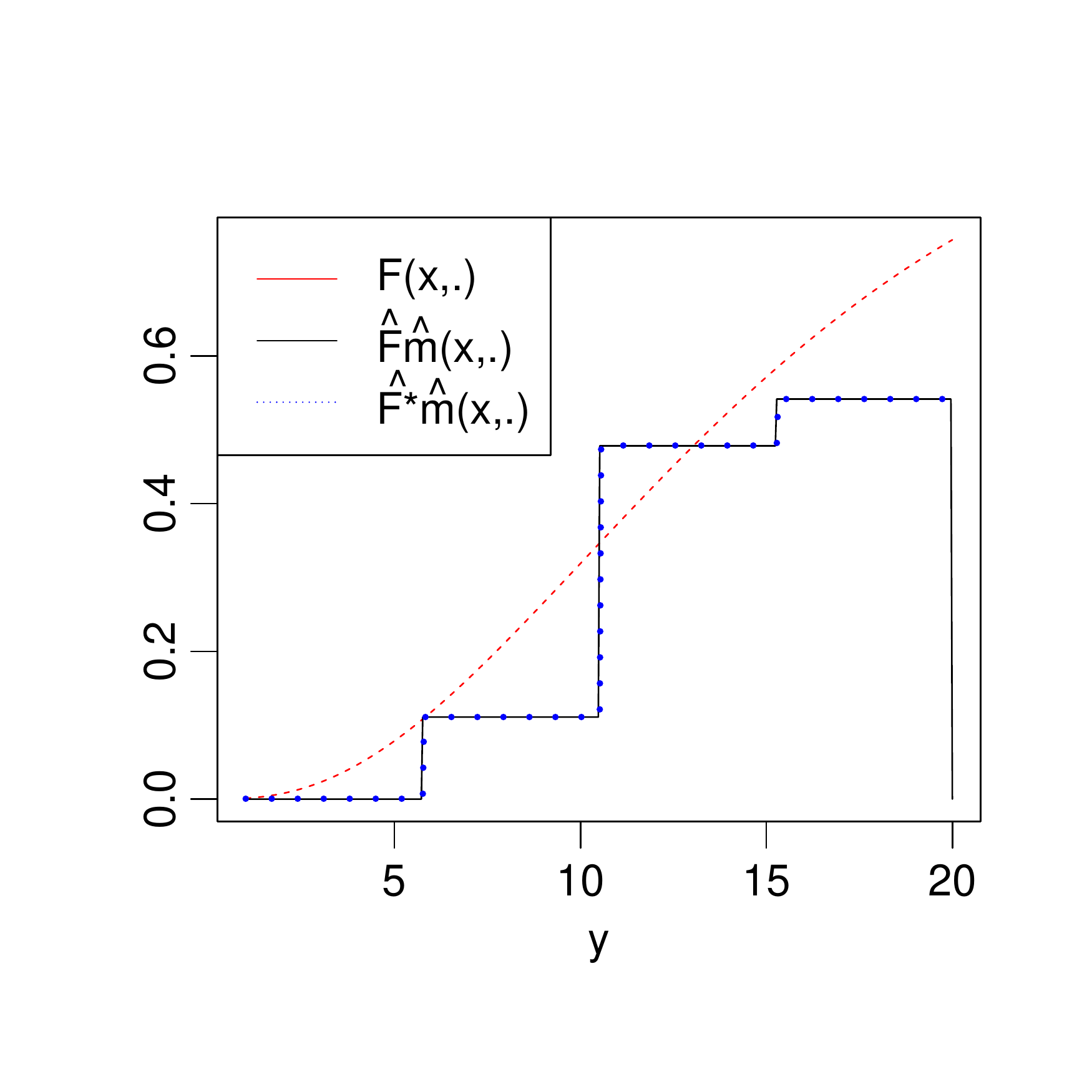}&
 \includegraphics[width=5cm,trim= 20mm 20mm 10mm 20mm]{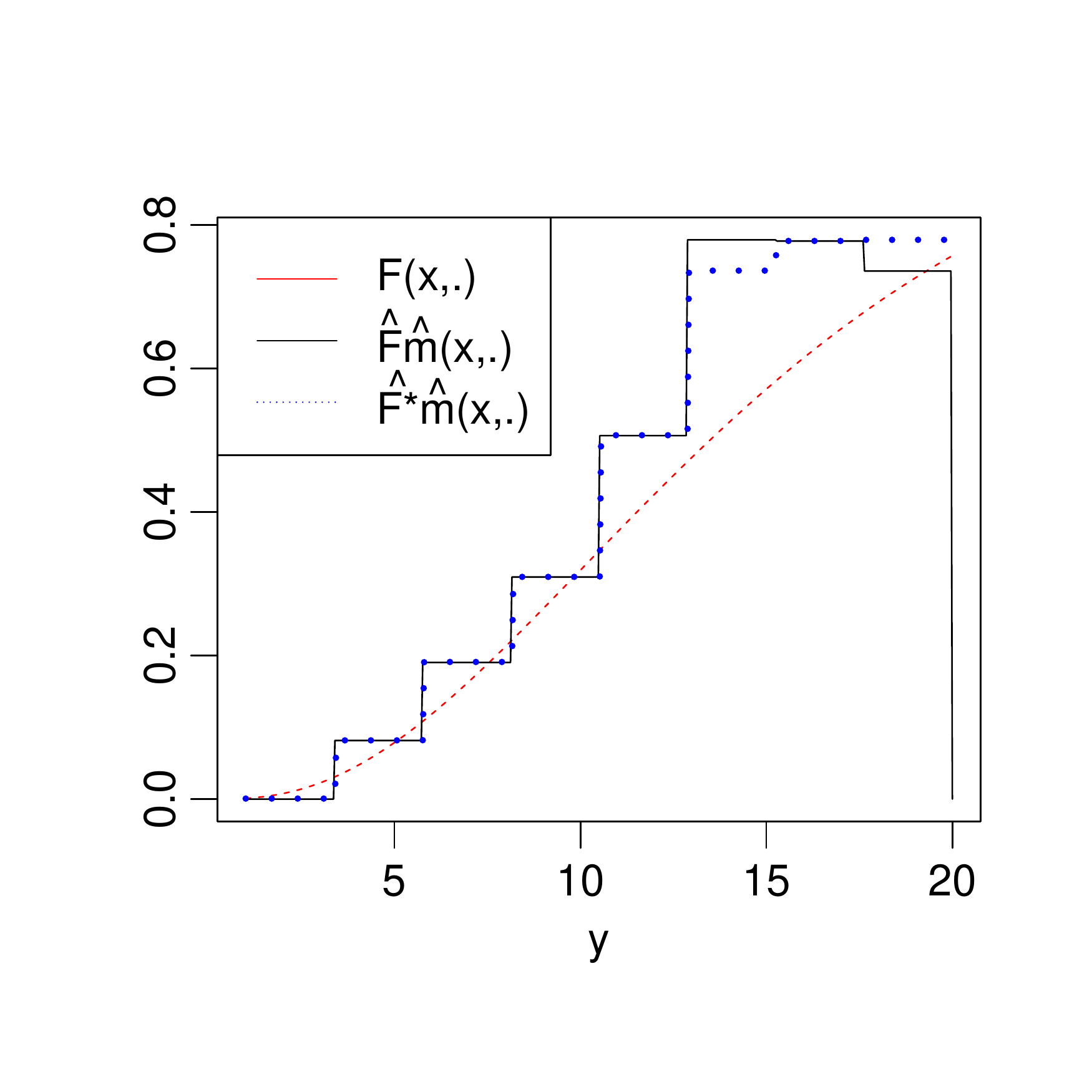} \\
 $n=500$ & $n=1000$ & $ n=5000$ 
\end{tabular}
\caption{\label{Fig-rea} Model selection estimator $\tilde{F}_{\widehat{m}}$ (black solid line), rearranged estimator $\tilde{F}_{\widehat{m}}^*$ (blue dotted line) and target function $F$ (red dashed line) for a fixed value of $x$. First row: model 1, $x=0.9$; second row: model 2, $x=5$.}
\end{figure}

The results displayed on Figure \ref{Fig-plotx} enhance that the estimate of $F$ is not necessarily increasing with respect to the second variable. Indeed the mean-square estimation does not account for the monotonicity of the c.d.f. The quality of estimation can be improved by monotonizing the model selection estimator $\tilde{F} _{\widehat{m}}$. For an univariate function $h$ defined on $[0,1]$, \cite{Cherno07} propose a rearrangement based on the following observation: if $h$ was non-decreasing, the quantile function of the random variable $h(U)$ with $U$ uniformly distributed on $[0,1]$ would be equal to $h$. Therefore, the rearrangement $h^*$ of $h$ is defined as the quantile function of $h(U)$:

\begin{equation*}
h^*(y) = \inf \left\{ z \in \mathbb{R}, \left( \int _0^1 1\!\text{I} _{\{ h(u)\leq z \} } du \right) \geq y \right\}.
\end{equation*}
For each $x \in A_1$, we implement this procedure to the marginal function $y \rightarrow \tilde{F}_{\widehat{m}} (x,y)$  up to an affine substitution of variable which transform $A_2=[a,b]$ into $[0,1]$. Namely,

\begin{equation}\label{eq-rea}
 \tilde{F}_{\widehat{m}} ^*(x,y) =  \inf \left\{ z \in \mathbb{R}, \left( \int _a^b 1\!\text{I} _{\{ \tilde{F}_{\widehat{m}}(x,u) \leq z \} } du \right) \geq y -a \right\}, \quad \forall (x,y) \in A.
\end{equation}
In the numerical implementation, the infimum in (\ref{eq-rea}) is taken over the range of values of $\tilde{F}_{\widehat{m}} (x,\cdot)$. Moreover, as $\tilde{F}_{\widehat{m}} $ is piecewise constant, the rearrangement is identical on each interval and only a finite number of rearrangement has to be computed. The results are presented on Figure \ref{Fig-rea}.

By construction the function $\tilde{F}_{\widehat{m}} ^*$ is increasing with respect to the second variable. Moreover, according to Proposition 1 in \cite{Cherno07}, the rearranged estimate presents a smaller $L^2$-error than the initial estimate: for every $x \in A_1$,
\begin{equation}
\int _{A_2} \left( \tilde{F}_{\widehat{m}} ^*(x,y) - F(x,y)  \right)^2 dy \leq \int _{A_2} \left( \tilde{F}_{\widehat{m}} (x,y) - F(x,y)  \right)^2 dy \quad a.s.
\end{equation}
Thus the rearrangement brings an improvement in terms of $L^2$-risk.

\section{Discussion}\label{sec-comments}

\subsection{Summary of the results} In the interval censoring framework, case 1, \cite{comte09} observe that the c.d.f. is equal to the regression function of the observation times over the indicators and propose an estimator of the c.d.f. computed by minimization of a least-square contrast on a collection of linear models followed by a model selection procedure. In this paper, we extend this approach to conditional c.d.f. via a collection of two-dimensional linear models built as tensor products of uni-dimensional ones. The inclusion of a covariate does not give rise to specific technical difficulties, but in order to emphasize the role played by each hypothesis, we adopt a different presentation of the results which leads to changes in the proofs. 


First of all, an oracle inequality for the empirical risk conditional to the observations $({\bf X,T})$ is stated. This result only requires a limitation of the number of models in the collection with a given dimension and is valid regardless of the regularity of the variables. Besides, it directly handles the fix-designed context with non-random observation times and covariates. The main difference in the proof is the use of a concentration inequality for empirical process with independent but non identically distributed variables.  An oracle inequality for the $L^2$-risk on a restricted set $A\in \mathbb{R}\times \mathbb{R}^+$ is then derived under additional assumptions. First of all, the density $f_{(X,T)}$ is assumed to be bounded on $A$ in order to guarantee the equivalence between the empirical and $L^2$-norms; this condition is inherent to the design since the c.d.f can only be estimated on a set with sufficiently concentrated observations. Moreover, the construction of our collection of models by tensor product requires the set $A$ to be  a product of intervals. 

The oracle inequalities state the adaptivity of the model selection estimator over the collection of models.  In order to demonstrate the optimality in a more general sense, a minimax study is conducted over anisotropic Besov balls. The rate is characteristic of regression estimation but differs from the rate obtained on right-censored data. 

A numerical study on simulated data with additive and multiplicative effects of the 
covariate is presented. We observe that a large sample size is necessary in order to get acceptable quality of estimation, due on the one hand to the multi-dimensional context, and on the other hand to the current status data framework in which poor information is brought by the censoring indicator.  Besides, in alternative to the censoring rate usually considered in right-censoring, we emphasize the impact of the distance between the densities of the survival and observation times on the quality of estimation: closer densities lead to a better estimation.

\subsection{Extension to multi-dimensional covariates}

From a theoretical point of view, the procedure presented in this paper could be extended to covariates $X\in \mathbb{R}^k$ with $k\geq 2$ by building a collection of models from multiple tensor products. For an index $m = (m_1, \dots , m_{k+1})$, the following model would be considered:

$$ S_m = S_{m_1}^{(1)}\otimes \dots \otimes S^{(k+1)}_{m_{k+1}}$$
with $S_{m_1}^{(1)}, \dots,  S^{(k)}_{m_{k}}$ linear subspaces of $L^2(\mathbb{R})$ and $S^{(k+1)}_{m_{k+1}}$ linear subspace of $L^2(\mathbb{R}^+)$.\\

Nevertheless the simulation study with a uni-dimensional covariate emphasizes the necessity of a large sample, and additional covariates could only be considered in situations with a large sample available. 
Besides the control of the integrated risk requires restrictions on the maximum dimension of the models.
Let $N_n^{(\ell)}$ be the maximum dimension of $S^{(\ell)}_{m_\ell}$ for $\ell = 1, \dots, k+1$, assumption $(\textbf{A}_2)$ turns into:

\begin{equation}\label{gen-eq1}
 N_n^{(1)}\times \dots  \times N_n^{(k+1)} \leq \sqrt{\frac{n}{\log n}} .
 \end{equation}
which is satisfied as soon as:
$$ N_n^{(\ell)} \leq \left( \frac{n}{\log n}\right) ^{1/2(k+1)}, \quad \forall \ell =1 ,\dots, k+1.$$ 
Assume that $F$ is in $\mathcal{B} _{2,\infty}^{\beta}(A,L) $ the anisotropic Besov ball of regularity $\beta = (\beta _1, \dots, \beta _{k+1})$ on $L^2(\mathbb{R}^k \times \mathbb{R}^+)$. Under condition (\ref{gen-eq1}), similarly to (\ref{eq-trade}), the model selection estimator satisfies the following inequality. 

\begin{eqnarray}\label{gen-eq2}
&&\hspace{-2cm} \nonumber \mathbb{E} \left[ \| \tilde{F} _{\widehat{m}} -F   \|^2\right]   \\
 &&  \hspace{-2cm} \leq C \left\{  
 \left( D_{m_1}^{(1)} \right) ^{-2\beta _1} + \dots + 
\left( D_{m_{k+1}}^{(k+1)} \right) ^{-2\beta _{k+1}}  + 
\frac{D_{m_1}^{(1)} \dots D_{m_{k+1}}^{(k+1)} }{n}  
  \right\}  + \frac{C'}{n}
\end{eqnarray}
with $D_{m_{\ell}}^{(\ell)}$ the dimension of $S^{(\ell)}_{m_\ell}$. Thus, the model $\overline{m} = (\overline{m}_1, \dots, \overline{m}_{k+1})$ which realizes the trade-off in (\ref{gen-eq2})  satisfies:

\begin{equation}\label{gen-eq3}
D_{\overline{m}_{\ell}}^{(\ell)} \propto n^{1/(2\beta _\ell ( 1+ \sum _{j=1}^{k+1} (2\beta _j )^{-1} )}.
\end{equation}
 Conditions   (\ref{gen-eq1}) and  (\ref{gen-eq3}) are satisfied as soon as $\beta _\ell > (k+1)/2 $ for every $\ell = 1,\dots, k+1$, and the rate of convergence is  $n^{-2\overline{\beta}/(2\overline{\beta} + (k+1))}$ with $\overline{\beta}$ the harmonic mean. Therefore, the regularity conditions which guarantee the minimax adaptivity get stronger and the rate of convergence decreases as the dimension of the covariate increases.

\subsection{Contribution with respect to extant literature}

Most estimation procedures developed for current status data are based on the NPMLE originally proposed by \cite{groeneboom92} and the covariates are usually handled through semi-parametric models (e.g. \cite{Cheng11}).  Minimax results are assessed under the assumption of a continuous density for the survival and observation times (\cite{vandergeer93}). Nevertheless, very few adaptive procedures are proposed and they do not include covariates. In alternative to the NPMLE, \cite{ma} and \cite{comte09} enlighten the parallel between regression and current status data frameworks and propose several adaptive estimation procedures adapted from regression analysis, which exhibit rates of convergence typical of a regression context. 

In this paper, we develop the mean-square procedure from \cite{comte09}. Our first contribution is the inclusion of a covariate in a nonparametric adaptive estimation procedure. The second improvement is the computation of the minimax rate over anisotropic Besov balls which assesses the optimality of the regression-type approach in the current status setting; the minimax rate for independent survival and observation times can be derived as a particular case, proving the minimaxicity of the mean-square estimator from \cite{comte09}.

\subsection{Perspectives}\label{sec-comparison}

The parallel between current status and regression frameworks would allow the implementation of other adaptive methods. We briefly present some advantages and drawbacks of alternative approaches that may be considered. \\

 Lepski's method provides locally minimax results up to a logarithm loss unavoidable in pointwise adaptive  procedures, and more recent developments display oracle inequalities.  
 Adaptation over classes of regularity is achieved through a bias-variance compromise from kernel estimators (\cite{lepski1}). \cite{Gold08} propose a modification of the original method based on the ordering of the classes of regularity in order to estimate multi-dimensional function in a gaussian regression framework. Adaptation of this procedure for current status data framework is conceivable, but as for model selection, a specific study of the upper-bound would be necessary. Besides, the implementation of Lepski's procedure presents a high computational complexity, especially in mullti-dimensional context. Since the estimation of the conditional c.d.f. from current status data requires large sample size, more easily computable methods may be preferred. However, Lepski's method allows local adaptivity which may be favoured when the target function present irregular smoothness. \\

 Wavelet thresholding methods are widely used in practical situations, in particular in signal processing. The procedure consists in a decomposition in a multi-scale wavelet basis followed by a thresholding of the  
 smallest coefficients. The estimators reach the minimax rate of convergence over various classes of regularity, usually up to a logarithm loss (e.g. \cite{Kerk04}). \cite{massart} demonstrates the equivalence between wavelet thresholding and the {\it complete} variable selection procedure which falls into the model selection framework with a collection including a large number of models with the same dimension. Therefore, following the adaptation of tools developed for regression estimation by model selection to the current status framework, wavelet thresholding could be developed in future research. \\
 
 The mean-square contrast is naturally  adapted to the empirical norm and additional assumptions are required for the control of the $L^2$-risk. Conversely, quotient regression estimators presented e.g. in \cite{tsybakov} are directly adapted to the $L^2$ norm. \cite{comte09} show that the rate of convergence of the model selection quotient estimator from current status data depends on the minimum regularity of the c.d.f. and the observation time density. Therefore, this procedure does not reach the minimax rate if the c.d.f. is smoother than the observation time. \\
 
 As the conditional c.d.f. lies into $[0,1]$, procedures for bounded regression may be considered. In particular, we present a comparison of our method with the model selection  framework for bounded regression proposed by \cite{massart}. The general assumptions  are less restrictive than in the classical regression framework since the errors are not supposed to be independent from the design, but the procedure requires models which consist of bounded functions. 
 Two types of bounded models could be directly derived from the linear models $\{ S_m, m\in I_n \}$:
 
 \begin{itemize}
 \item Thresholded linear models: $\tilde{S}_m^{(1)} = \{ \tilde{t} = \max (\min (t,1),0), t\in S_m \}$
 
 \item Restricted linear models: $\tilde{S}_m^{(2)} =  \{ t\in S_m, t(x) \in [0,1], \forall x \}$. 
 
 \end{itemize}
We present a brief overview of the advantages and difficulties encountered in the implementation of these models on the interval censoring design. According to Theorem 8.5 in \cite{massart}, for each model $m$, two quantities have to be controlled  in order to obtain minimax results: 1/ the stochastic modulus of uniform continuity $\phi _m$ of the centered empirical process, 2/ the bias term via results of approximation on regularity spaces. On the one hand, the bias term for $\tilde{S}_m^{(1)}$ is directly controlled by results on linear models, but since the unit ball of $\tilde{S}_m^{(1)}$ is larger than the one of $S_m$, the computation of the function $\phi _m$ requires a specific study. On the other hand, as $\tilde{S}_m^{(2)}$ is included in $S_m$,  $\phi _m$ can directly be  deduced from results on linear models, but the bias term which involved the $L^2$-projection of $F$ on $L^{\infty}$-balls does not falls into classical approximation setting. 

Besides, an alternative procedure based on finite models, originally developed by \cite{massart} in a Gaussian framework (Chapter 4), could be considered in order to include bounded models.  For every $\alpha >0$, a net of radius $\alpha$ is extracted from a set of function $\mathcal{S}$ including the target function. An appropriate choice of radii $\{ \alpha _m , m\in I_n\}$ based on wavelet or polynomial characterization of Besov classes provides minimax estimation.

Thus, our model selection procedure takes advantage of classical tools for the control of the deviation of empirical process on linear models as well as extant results of approximation on Besov spaces, but does not directly fit into the bounded  regression framework.  Alternative approaches using bounded models could be considered for future research, but give rise to other technical difficulties than the one encountered with linear models.

\subsection{Conclusion} The model selection mean-square estimator of the conditional c.d.f. satisfies non-asymptotic adaptivity on a collection of models, as well as minimaxity. The control of the upper-bound takes advantages of the tools used in classical regression.  Furthermore, our results underline the optimality of the regression-type approach, which establishes a connection between the interval censoring framework and the widely developed field of regression estimation and makes the way for the adaptation of regression methods to current status design.

\section{Proofs}\label{sec-proofs}

\subsection{Talagrand Inequality}

The risk of the model selection estimator is controlled with the following concentration inequality based on results from \cite{talagrand}.

\begin{theorem}\label{talagrand-noniid}
Let $(V_1, \dots, V_n)$ be independent random variables, and $\mathcal{F}$ be a countable set of applications from $\mathbb{R}$ to $\mathbb{R}^n$ such that for every $f\in \mathcal{F}$, the function $-f$ defined as $(-f)(x)=-(f(x))$ for every $x$ is in $\mathcal{F}$ as well. Let 
$$Z= \sup _{f\in \mathcal{F}}\left|\frac{1}{n} \sum _{i=1}^n (f^{(i)} (V_i) - \mathbb{E} [ f^{(i)}(V_i)])\right| $$
and $b$, $v$ and $\mathbb{H}$ be such that

\begin{equation*}
\sup_{f\in \mathcal{F}} \left(\sup _{i=1,\dots ,n} \| f^{(i)} \|_{\infty} \right) \leq b, \quad \sup_{f\in \mathcal{F}} \frac{1}{n} \sum _{i=1}^n Var (f^{(i)}(V_i))\leq v \quad \text{and} \quad \mathbb{E}Z\leq \mathbb{H}.
\end{equation*} 
Then, for every $\theta >1$, there exists $\overline{C}$, $\overline{C}'$, $\overline{K}$, $\overline{K}'$ such that for every $n$,

$$\mathbb{E} \left[ (Z^2 - \theta \mathbb{H}^2)_+ \right] \leq \overline{C}\frac{v}{n}\exp \left(-\overline{\kappa} \frac{n\mathbb{H}^2}{v}\right) +\overline{C}'\frac{b^2}{n^2} \exp \left(-\overline{\kappa }' \frac{n\mathbb{H}}{b}\right)$$

\end{theorem}

Theorem $\ref{talagrand-noniid}$ is derived from Theorem 1.1 in \cite{klein-rio} by setting $s^{(i)}(t) = \frac{1}{b} ( f^{(i)} (t) - \mathbb{E} [ f^{(i)}(V_i)] )$. Similarly to \cite{BM}, Corollary 2, we get
$$P[ Z \geq (1+ \nu) \mathbb{H} +x] \leq \exp \left( -\frac{n}{3} \min \left( \frac{y^2}{2v} , \frac{ 2 \min (1,\nu) y}{7b} \right) \right).$$ 
Then we use the formula $\mathbb{E} [X_+ ] \leq \int _0^{+\infty} P[ X \geq s] ds$ to get
$$\mathbb{E} \left[ (Z^2 - \theta \mathbb{H}^2)_+ \right] \leq \int _0^{+\infty} \!\!\! P[Z \geq \sqrt{ \theta \mathbb{H} ^2 +s } ] ds \leq \int _0^{+\infty} \!\!\!P[Z \geq \alpha _1 \mathbb{H} + \sqrt{ \alpha _2 \mathbb{H} + \alpha _3 s}  ] ds $$
for some positive $(\alpha_1 , \alpha _2, \alpha _3)$, and a simple integration provides the result of Theorem \ref{talagrand-noniid}.

\subsection{Proof of Theorem $\ref{T1}$. }

The proof is similar to \cite{comte09}, Section 6.4, but the conditional empirical norm is considered instead of the $L^2$ one. This modification entails changes in the proofs that are displayed in this section. The vectors ${\bf X}$  and ${\bf T}$ are assumed to be fixed along this section. Let $m \in I_n$ and $F_m \in S_m$, by definition of $\widehat{F}_m$ and  $\widehat{m}$,

\begin{eqnarray*}
 \| \widehat{F}_m - F \|_n ^2 \hspace{-0.2cm}& \leq  \hspace{-0.2cm}&  \|F _m - F\|_n ^2 + pen (m) - pen (\widehat{m}) + 2 \nu _n (\widehat{F} _{\widehat{m}} - F_m ) \\
 \hspace{-0.2cm}& \leq  \hspace{-0.2cm}&\|F _m - F\|_n ^2 + pen (m) - pen (\widehat{m}) + \| \widehat{F}_m - F_m  \|_n \sup _{t \in S_m + S_{\widehat{m}}, \|t\|_n\leq 1}\!\!\! | \nu _n (t)|
\end{eqnarray*}
with $\nu _n$ defined in (\ref{eq-nu}).  Let $\theta ' = \theta ^{1/3} >1$, and $p(m,m')= \theta ' (D_m + D_{m'})/(4n)$. We recall that for every $(x,y) \in \mathbb{R}^2$, $2xy \leq x^2/\theta'^2 + \theta '^2 y^2$ and $(x+y)^2 \leq \theta ' x^2 + (1+1/(\theta' -1)) y^2$. Therefore,

\begin{eqnarray*}
 \| \widehat{F}_m - F \|_n ^2 
&\leq & \|F _m - F\|_n ^2 + pen (m) - pen (\widehat{m}) +  \theta '^2 p(m, \widehat{m})  + \frac{1}{\theta '^2} \| \widehat{F}_m - F_m  \|_n ^2  \\
&& + \theta '^2 \sup _{t \in S_m + S_{\widehat{m}}, \|t\|_n\leq 1} \left( (\nu _n (t))^2 - p(m, \widehat{m}) \right) \\
&    \leq & \|F _m - F\|_n ^2 + 2pen (m)  \\
&& +  \frac{1}{\theta '^2} \left( \| \widehat{F}_m - F  \|_n ^2 
+ \theta ' \left( 1+ \frac{1}{\theta '-1} \right) \| F_m - F  \|_n ^2   \right) \\
&&   + \theta '^2  \sum _{m' \in I_n} \sup _{t \in S_m + S_{m'}, \|t\|_n\leq 1} \left( (\nu _n (t))^2 - p(m, m') \right)
\end{eqnarray*}
The end of the proof is provided by Theorem \ref{talagrand-noniid}. According to the comment in Section \ref{sec-leastsquare}, the supremum of $\nu _n $ over the $\|.\|_n$-unit ball of $S_m + S_{m'}$ is equal to the supremum over the unit ball $\mathcal{B}_{m+m'}^*$ of $(S_m + S_{m'})^*$. Let $\{ \varphi _\lambda , \lambda \in J_{m+m'}^* \}$ be a $\|.\|_n$-orthogonal basis of $(S_m + S_{m'})^*$, by Cauchy-Schwartz Inequality and equation (\ref{eq-bias}),

\begin{equation}
\mathbb{E} \left[   \sup _{t  \in \mathcal{B}_{m+m'}^*} (\nu _n(t))^2 | ({\bf X,T})  \right] \leq \mathbb{E} \left[   \sup _{\lambda \in J_{m+m'}^*} (\nu _n(\varphi _\lambda ))^2 | ({\bf X,T})  \right] \leq \frac{D_m + D_{m'}}{4n} = \mathbb{H}^2. 
\end{equation}
Moreover, for every $t  \in \mathcal{B}_{m+m'}^*$ and $ i \in \{ 1,\dots, n\}$,
\begin{equation}
| \delta _i t(X_i, T_i) | \leq | t(X_i,T_i)| \leq \sqrt{n}  \| t\| _n  \leq \sqrt{n} = b. 
\end{equation}
Finally, with a similar computation as (\ref{eq-bias}), the supremum of the variance process is upper bounded by $v=1/4$. Therefore, Theorem \ref{talagrand-noniid} with the values of $\mathbb{H}$, $b$ and $v$ computed above concludes the proof of Theorem \ref{T1}. $\Box$

\subsection{Proof of Theorem \ref{T2-cens-int}}

The proof is based on results presented by \cite{baraud-random} in a uni-dimensional context, but the demonstrations are identical for two-dimensional functions. Let $\rho _0 > h_0^{-1}$ and 

$$ \Omega _n = \left\{ \frac{\|t\|_n^2}{ \| t\|^2} \leq \rho _0, \ \forall t \in S_n  \right\} . $$
On the one hand, according to Lemma 5.1 in \cite{baraud-random}, 

\begin{equation}\label{eq1-T2}
\mathbb{E} \left[  \| \tilde{F} _{\widehat{m}} -F \| 1\!\text{I} _{\Omega _n} \right] \leq
C \inf _{t \in S_n} \| F-t \| ^2 + \mathbb{E} \left[  \| \tilde{F} _{\widehat{m}} -F \|_n ^2  \right]
\end{equation}
for some constant $C$. Moreoever, by integrating (\ref{eq-T1-chap5}) in Theorem \ref{T1}, we obtain:

\begin{equation}\label{eq2-T2}
 \mathbb{E} \left[  \| \tilde{F} _{\widehat{m}} -F \|_n ^2  \right] \leq C_1 \inf _{ m \in I_n} \left\{ \inf _{t \in S_m} h_1 \| F-t \| ^2 + pen (m) \right\} + \frac{C_2}{n}.
 \end{equation}
Thus, as $\inf _{t \in S_n} \| F-t \| ^2 \leq \inf _{t \in S_m} \| F-t \| ^2 $ for all $m \in I_n$, inequalities (\ref{eq1-T2}) and (\ref{eq2-T2}) entail:

\begin{equation}\label{eq3-T2}
\mathbb{E} \left[  \| \tilde{F} _{\widehat{m}} -F \| 1\!\text{I} _{\Omega _n} \right] \leq C' \inf _{ m \in I_n} \left\{ \inf _{t \in S_m}  \| F-t \| ^2 + pen (m) \right\} + \frac{C_2}{n}.
\end{equation}
On the other hand, according to Proposition 4.2 in \cite{baraud-random},

\begin{equation}\label{eq4-T2}
\mathbb{P} \left[ \Omega _n ^c\right] \leq \frac{C''}{n}.
\end{equation}
Moreover, $\tilde{F} _{\widehat{m}} $ and $F$ lie into $[0,1]$, therefore

\begin{equation}\label{eq5-T2}
 \| \tilde{F} _{\widehat{m}} -F \| \leq h_0^{-1} \| \tilde{F} _{\widehat{m}} -F \|_{f_{X,T}} \leq  h_0^{-1} \quad a.s.
\end{equation}
Equations (\ref{eq3-T2}), (\ref{eq4-T2}) and (\ref{eq5-T2}) conclude the proof of Theorem \ref{T2-cens-int}.  $\Box$


\subsection{Proof of Theorem \ref{T-low}}
 The proof is based on the following theorem and lemma presented in \cite{tsybakov} (Theorem 2.5 and   Lemma 2.7 in Chapter 2). 
  Let $\mathcal{B}= \mathcal{B}^{\beta}_{2, \infty} (A,L)$. Denote by $K(P,Q)$ the Kullback distance between the distributions $P$ and $Q$:

 $$ K(P,Q)= 
 \left\{ \begin{array}{l}
 \int \log (dP/dQ) dP \quad  \text{if}  \quad P << Q \\
 +\infty \quad \text{otherwise} 
\end{array}\right. $$

\begin{theorem}\label{th-tsy}

Assume that there exist $M\geq 2$ and $F_0, \dots , F_M$ such that

\begin{enumerate}

\item $F_j \in \mathcal{B}$ for every $j\in \{ 0, \dots , M\}$.

\item $\| F_j -F_l \| ^2 \geq 2r$ for every $j \neq l \in \{ 0, \dots , M\} .$

\item $P_j^{(n)} << P_0 ^{(n)}$ for every $j\in \{ 0, \dots , M\}$, where $P_j^{(n)} $ denotes the distribution of $(X_i,T_i,\delta _i )_{i=1, \dots, n}$ if $F= F_j$, and for some  $0 < \alpha < 1/8$
$$\frac{1}{M} \sum _{j=1}^M K( P_j^{(n)} , P_0^{(n)} ) \leq \alpha \log M.$$
\end{enumerate}
Then there exists a constant $c$ such that
$$\inf_{\widehat{F}_n } \sup _{F\in \mathcal{B} } \mathbb{E} \left[ r \| \widehat{F}_n - F \| ^2 \right] \geq c.$$

\end{theorem}

\begin{Lemma}\label{Lemma-Tsy}
Let $m$ be a positive integer. There exists a family $(\varepsilon ^{(0)}, \dots, \varepsilon ^{(M)} ) \subset \{0,1 \} ^{2^m} $ with $\varepsilon ^{(0)}= (0, \dots, 0)$ such that $\log (M) \geq (\log 2/8) 2^m$, and
$$\overline{\rho}(\varepsilon ^{(i)} , \varepsilon ^{(i')}) \geq  \frac{2^m}{8}, \ \ \forall i\neq i' \in \{0, \dots ,M\}$$
with  $\overline{\rho}$ denoting the following distance:
\begin{equation}\label{def-rho-chap6}
 \overline{\rho} (\varepsilon , \varepsilon ') =  \sum _{k =1}^{2^m}\text{1\!{\mbox I}}_{\{ \varepsilon _k \neq \varepsilon '_k\}}.
\end{equation}

\end{Lemma}

Let us present the outline of the proof. Up to rescaling and translation, we assume that $A=[0,1]\times [0,1]$. We build a set of c.d.f.:
$$ F_{\varepsilon} = F_0 + \sqrt{\frac{b}{n}} \sum _{S\in R_J} \varepsilon _S \psi _{J,S} $$
where $F_0$ is a baseline c.d.f.,   $\varepsilon_S \in \{0,1\}$ and the $\{ \psi _{J,S}, S\in R_J\} $ are wavelets with disconnected supports included in $[0,1]$, defined from a mother wavelet $\psi$.

\begin{equation}\label{eq-wave}
\psi_{J,S}(x,u)= 2^{(j_1+j_2)/2} \psi (2^{j_1} x-s_1) \psi (2^{j_2} u -s_2)
\end{equation}
for some $J=(j_1,j_2)\in \mathbb{N}^2$ and $R_J \subset \mathbb{Z}^2$. We show that

\begin{itemize}
\item  $ \| F_{\varepsilon} - F_{\varepsilon '} \|^2=  (b/n)  \overline{\rho} (\varepsilon , \varepsilon ').$

\item $ K( P_{\varepsilon} ^{(n)} , P_0 ^{(n)}) \leq b'  2^{j_1 + j_2} $.

\end{itemize}
Therefore, for every $J$, Lemma \ref{Lemma-Tsy} ensures the existence of a set $( F_{\varepsilon ^{(0)}}, \dots, F_{\varepsilon ^{(M)}})$ such that
conditions $(2)$ and $(3)$ in Theorem \ref{th-tsy} are satisfied with an appropriate choice of $b$. 
Moreover, we show that $\| F_{\varepsilon}  \| _{\mathcal{B}^{\beta} _{2, \infty} ( [0,1]^2)}$ has order 
\begin{equation}\label{order}
\frac{2^{(j_1+j_2)/2} (2^{j_1\beta _1} + 2^{j_2 \beta _2}+1) }{\sqrt{n}}
\end{equation}
 Then $J$ is chosen is order to guarantee that (\ref{order}) is upper bounded by a constant, which leads to the rate of convergence $r = 2^{j_1+j_2} \propto n^{- \overline{\beta} /(\overline{\beta} + 1)}$. We now present the detailed proof.

\subsubsection{Construction of the $(F_i)$'s.}

Let
$$F_0(x,u)=1\!{\mbox I}_{[0,1]} (x) \left( a1\!{\mbox I}_{[0, +\infty[ }(u) + a u1\!{\mbox I}_{[0,1]}(u) + (1-a)1\!{\mbox I}_{(1, + \infty[} (u) \right)$$
with $a= \min (1/3, L/2)$.
For every $x\in [0,1]$,
\begin{eqnarray*}
&& \bullet \  F_0(x,u)= 0, \ \ \forall u <0,\\
&& \bullet \   F_0(x,u)=1, \ \ \forall u > 1, \\
&& \bullet \  F_0(x,.) \ \text{ is increasing on  }[0,1] \text{  and  }F_0(x,u) \in [a,2a] \subset(0,1) \text{  for every  } u\in [0,1],
\end{eqnarray*}
thus $F_0$ is a conditional distribution. Let $\psi$ be a one-dimensional wavelet supported on $[0,1]$. Let $J=(j_1,j_2) $ be a couple of non-negative integers determined further. For every $S=(s_1, s_2) \in \mathbb{Z}^2$, let $psi_{J,S} $ be defined in (\ref{eq-wave}).
There exists a subset $R_J$ of $\mathbb{Z}^2$ such that
\begin{eqnarray*}
&& \bullet \  Supp(\psi _{J,S})= I_{J,S} \subset ]0,1[^2 \text{  for every  }S\in R_J,\\
&& \bullet \  \text{The applications  }\{ \psi _{J,S} , S\in R_J \} \text{  have disjoint supports }, \\
&& \bullet \  |R_J|= 2^{j_1+j_2}.
\end{eqnarray*}
Let $b$ be a positive constant which will be determined later.
For every $\varepsilon \in \{ 0,1 \} ^{| R_J |}$, let
$$G_{\varepsilon} = \sqrt{\frac{b}{n}} \sum _{S\in R_J} \varepsilon _S \psi _{J,S} $$
and $F_{\varepsilon} = F_0 + G_{\varepsilon}$. For every
$x\in [0,1]$,
\begin{eqnarray*}
&& \bullet \  F_{\varepsilon}(x,u)= F_0(x,u)=0, \ \ \forall u <0 \\
&& \bullet \  F_{\varepsilon}(x,u)= F_0(x,u)=1\ \ \forall u > 1
\end{eqnarray*}
Moreover let $(x,u)\in[0,1]^2 $,
\begin{equation}\label{int-F}
F_{\varepsilon} (x,u) =  a + \int _0^u \left(a + \sqrt{\frac{b}{n} } \sum _{S\in R_J} \varepsilon _S \frac{\partial \psi _{J,S} }{\partial y} (x,y) \right) dy.
\end{equation}
Assume that
\begin{equation}\label{C1}
\sqrt{\frac{b}{n}} 2 ^{j_1/2} 2^{3j_2/2} \| \psi \| _{\infty} \left\| \psi ' \right\| _{\infty} \leq \frac{a}{2}. 
\end{equation}
Let $y \in [0,u]$ and $S_0$ be such that $(x,y)\in I_{J,S_0}$
$$\left|  \sqrt{\frac{b}{n} } \sum _{S\in R_J} \varepsilon _S \frac{\partial \psi _{J,S} }{\partial y} (x,y) \right| = \left| \sqrt{\frac{b}{n} } \varepsilon _{S_0} \frac{\partial \psi _{J,S_0} }{\partial y} (x,y) \right| < \frac{a}{2}.$$
 Therefore the term in the integral in ($\ref{int-F}$) is positive and the application $F_{\varepsilon} (x,.) $ is increasing on $[0,1]$. Moreover, as $\psi _{J,S}(x,1)=0$ for every $S\in R_J$, $ F_{\varepsilon} (x,1) = F_0(x,1)  = 2a <1$.
Thus $F_{\varepsilon}$ is a conditional distribution function on $[0,1]^2$.

\subsubsection{Condition which guarantees that $F_{\varepsilon} \in\mathcal{B}$ for every $\varepsilon$.}

On the one hand, assume that $\psi$ is regular enough, then according to \cite{hochmuth} (Theorem 3.5),
\begin{eqnarray*}
| G_{\varepsilon} | _{\mathcal{B}^{\beta} _{2, \infty} ( [0,1]^2)} \leq (2^{j_1\beta _1} + 2^{j_2 \beta _2} ) \| G_{\varepsilon} \|.
\end{eqnarray*}
Moreover, as the $\{ \psi _{J,S} , S\in R_J\}$ have disjoint supports,
$$\| G_{\varepsilon} \| ^2 = \frac{b}{n} \left\|  \sum _{S\in R_J} \varepsilon _S \psi _{J,S} \right\|^2 = \frac{b}{n}  \sum _{S\in R_J} \varepsilon _S ^2 \left\| \psi _{J,S} \right\|^2.$$
By definition of the wavelets, $\|\psi _{J,S} \| = \| \psi \| = 1$, hence
$$\| G_{\varepsilon} \| \leq \sqrt{\frac{b}{n}}  | R_J| =  \sqrt{\frac{b}{n}} 2^{(j_1+j_2)/2}.$$
Thus
\begin{eqnarray*}
&&\| G_{\varepsilon} \| _{\mathcal{B}^{\beta} _{2, \infty} ( [0,1]^2)}  = | G_{\varepsilon} | _{\mathcal{B}^{\beta} _{2, \infty} ( [0,1]^2)} + \| G_{\varepsilon} \|   \leq  \sqrt{\frac{b}{n}} 2^{(j_1+j_2)/2} (2^{j_1\beta _1} + 2^{j_2 \beta _2}+1 ).
\end{eqnarray*}
On the other hand, $| F_0  | _{\mathcal{B}^{\beta} _{2, \infty} ( [0,1]^2)} =0$. Indeed, let $r_i = \lfloor \beta _i \rfloor +1$ for $i=1$ and $2$. Then $r_1 \geq 1$, $r_2 \geq 2$ and

$$| F_0  | _{\mathcal{B}^{\beta} _{2, \infty} ( [0,1]^2)} = \sup _{t>0} \left[ t^{-\beta _1} \omega _{r_1,1} (F_0,t,[0,1]^2)_2 + t^{-\beta _2} \omega _{r_2,2} (F_0,t,[0,1]^2)_2 \right].$$
Besides let $h>0$ and 
$$\Omega _{h,1}^{r_1} = \{ (x,u) \in [0,1]^2, (x+r_1h,u)\in [0,1]^2 \}.$$
For every $(x,u) \in \Omega _{h,1}^{r_1}$ $F_0(x+h,u)= F_0(x,u)$. So, as $r_1\geq 1$, $\Delta _{h,1}^{r_1} F_0 (x,u) = 0.$
Hence
$$ \omega _{r_1,1} (F_0,t,[0,1]^2)_2 = \sup _{|h|\leq t } \|  \Delta _{h,1}^{r_1} F_0 \| _{L^2(\Omega _{h,1}^{r_1})} =0.$$
Moreover, for every $(x,u) \in [0,1]^2 $
\begin{equation*}
\Delta _{h,2}^1 F _0 (x,u)= ah \quad \Rightarrow  \quad \Delta _{h,2}^{r_2} F _0 (x,u) = \Delta _{h,2}^{r_2-1}\Delta _{h,2}^1 F _0 (x,u) =0
\end{equation*}
as $r_2-1\geq 1$. Then $ \omega _{r_2,2} (F_0,t,[0,1]^2)_2=0$ and consequently $| F_0  | _{\mathcal{B}^{\beta} _{2, \infty} ( [0,1]^2)} =0$. Moreover

$$\| F_0  \| _{\mathcal{B}^{\beta} _{2, \infty} ( [0,1]^2)} = \sqrt{ \int _0^1 \int _0^1 a^2 (1+u)^2 du dx } = \sqrt{\frac{7}{3}} a $$
and
$$\| F_{\varepsilon}  \| _{\mathcal{B}^{\beta} _{2, \infty} ( [0,1]^2)} \leq  \| F_0 \| _{\mathcal{B}^{\beta} _{2, \infty} ( [0,1]^2)}+\| G_{\varepsilon}  \| _{\mathcal{B}^{\beta} _{2, \infty} ( [0,1]^2)} \leq \sqrt{\frac{7}{3}} a +  \sqrt{\frac{b}{n}} 2^{(j_1+j_2)/2} (2^{j_1\beta _1} + 2^{j_2 \beta _2}+1 ).$$
By definition $a \leq L/2$ so $\| F_{\varepsilon}  \| _{\mathcal{B}^{\beta} _{2, \infty} ( [0,1]^2)} \leq L$ as soon as
\begin{equation}\label{C2}
\sqrt{\frac{b}{n}} 2^{(j_1+j_2)/2} (2^{j_1\beta _1} + 2^{j_2 \beta _2}+1 ) \leq L \left(1- \sqrt{ \frac{7}{12}} \right).
\end{equation}

\subsubsection{Expression of $\| F_{\varepsilon} - F_{\varepsilon '} \| ^2$.}

\begin{equation}\label{def-rho-chap6}
 \| F_{\varepsilon} - F_{\varepsilon '} \|^2=  \frac{b}{n} \sum _{S\in R_J}  \int _{I_{J,S}} (\varepsilon _S - \varepsilon '_S) ^2 \psi _{J,S}^2 (x,u) dx du  =  \frac{b}{n} \sum _{S\in R_J}1\!{\mbox I}_{\{ \varepsilon _S \neq \varepsilon '_S\}} = \frac{b}{n} \overline{\rho} (\varepsilon , \varepsilon ').
\end{equation}

\subsubsection{Upper bound of $K(P_{\varepsilon } ^{(n)} , P_0^{(n)} )$}
For every $i \in \{1, \dots ,n\}$, under $F_{\varepsilon}$, $(X_i, T_i, \delta _i)$ has density

$$p_{\varepsilon}(x,u,d)= \left[ (F_{\varepsilon} (x,u) )^d (1- F_{\varepsilon} (x,u))^{1-d} \right]  f_{(X,T)} (x,u)$$
with respect to $\mathcal{L} \otimes \mathcal{L} \otimes \mu$ where $\mathcal{L}$ is the Lebesgue measure and $\mu$ is the counting measure on $\mathbb{N}$. Similarly, under $F_0$, $(X_i, T_i, \delta _i)$ has density

$$p_0(x,u,d)= \left[ (F_0 (x,u) )^d (1- F_0(x,u))^{1-d} \right]  f_{(X,T)} (x,u)$$
with respect to $\mathcal{L} \otimes \mathcal{L} \otimes \mu$. For every $\varepsilon \in \{ 0,1 \} ^{|R_J|}$, $P_{\varepsilon}$ is absolutely continuous with respect to $P_0$. Indeed,
$$F_0(x,u)=0 \quad \Rightarrow \quad  (x,u) \notin [0,1] \times [0, +\infty[ \quad  \Rightarrow \quad  F_{\varepsilon}(x,u)=0,$$
$$ F_0(x,u)=1  \quad \Rightarrow \quad  (x,u) \in [0,1] \times [1, +\infty[ \quad  \Rightarrow \quad  F_{\varepsilon}(x,u)=1.$$
Then
\begin{eqnarray*}
&&\hspace{-0.7cm} K(P_{\varepsilon} , P_0) =\\
 & & \hspace{-0.7cm} \int _{\mathbb{R}^2} \left[ \log \left( \frac{ F_{\varepsilon} (x,u)}{ F_0(x,u)} \right) F_{\varepsilon} (x,u) + \log  \left( \frac{ 1-F_{\varepsilon} (x,u)}{1- F_0(x,u)} \right)(1- F_{\varepsilon} (x,u)) \right] f_{(X,T)} (x,u) dx du
\end{eqnarray*}
Out of the intervals $\{ I_{J,S}, S \in R_J\}$, $F_{\varepsilon}$ and $F_0$ are equal. Hence
\begin{eqnarray*}
K(P_{\varepsilon}, P_0) & = & \sum _{S\in R_J}  \int _{I_{J,S}} \left[  \log\left( 1+ \frac{ \theta _S}{a(1+u)} \right) ( a(1+u) + \theta _S) \right.\\
&& \quad  + \left. \log\left( 1- \frac{ \theta _S}{1-a(1+u)} \right) (1- a(1+u) - \theta _S)\right] f_{(X,T)} (x,u) dx du
\end{eqnarray*}
where $\theta _S = \varepsilon _S \sqrt{ b/n } \ \psi _{J,S} (x,u).$
For every $S\in R_J$ and $(x,u) \in I_{J,S}$
$$ \frac{ \theta _S}{a(1+u)}  = \frac{ F_{\varepsilon} (x,u)}{ F_0(x,u)}-1 >-1 \quad \text{and} \quad - \frac{ \theta _S}{1-a(1+u)} = \frac{ 1-F_{\varepsilon} (x,u)}{1- F_0(x,u)}-1 >-1$$
Noting that $\log(1+v) \leq v$ for every $v>-1$ we obtain

$$K(P_{\varepsilon} , P_0) \leq \sum _{S\in R_J}  \int _{I_{J,S}} \left[ \theta _S + \frac{ \theta _S^2}{a(1+u)} - \theta _S + \frac{ \theta _S ^2}{ 1-a(1+u)} \right] f_{(X,T)} (x,u) dx du .$$
For every $u\in [0,1]$,
$$0 <\frac{1}{a(1+u)} \leq \frac{1}{a} \quad \text{ and } \quad 0 < \frac{ 1}{ 1-a(1+u)} \leq \frac{1}{1-2a}.$$
Thus
\begin{eqnarray*}
K(P_{\varepsilon} , P_0) &\leq& \left(\frac{1}{a} +\frac{1}{1-2a}\right) \sum _{S\in R_J}  \int _{I_{J,S}} \theta _S^2 f_{(X,T)} (x,u) dx du \\
& \leq  & \left(\frac{1}{a} +\frac{1}{1-2a}\right) \frac{b}{n} \| f_{(X,T)}\| _{\infty} |R_J| 
=   a 'b\| f_{(X,T)}\| _{\infty} \frac{ 2^{j_1+j_2}}{n} .
\end{eqnarray*}
where $a'= 1/a+ 1/(1-2a)$. Finally,

\begin{equation*}
K( P_{\varepsilon} ^{(n)} , P_0 ^{(n)}) \leq a'b \| f_{(X,T)}\| _{\infty}  2^{j_1 + j_2} .
\end{equation*}

\subsubsection{Conclusion}

According to Lemma \ref{Lemma-Tsy}, there exists a family $(\varepsilon ^{(0)}, \dots, \varepsilon ^{(M)} ) \subset \{0,1 \} ^{|R_J|} $ with $\varepsilon ^{(0)}= (0, \dots, 0)$ such that $\log (M) \geq (\log 2/8) 2^{j_1+j_2}$ and
$$\overline{\rho}(\varepsilon ^{(i)} , \varepsilon ^{(i')}) \geq \frac{|R_J|}{8} = \frac{2^{j_1+j_2}}{8}, \ \ \forall i\neq i' \in \{0, \dots ,M\}.$$  

Now the parameters $B_0$, $b$, $j_1$ and $j_2$ are choosen so that the family $(F_{\varepsilon ^{(0)}}, \dots, F_{\varepsilon ^{(M)}}) $ satisfies the assumptions of Theorem $\ref{th-tsy}$ with
$$r =B_0 n^{ \overline{\beta}/(\overline{\beta} +1)}.$$
Let
$$b= \frac{ \log 2}{72 \| f_{(X,T)} \| _{\infty} a'}  , \quad  \quad c_0= \left[ \frac{L}{4\sqrt{b}}\left(1- \sqrt{\frac{7}{12}}\right) \right]^{1/(1+ \beta _1 + \beta _2)}  \hspace{0.2cm} \text{and} \quad B_0=b c_0^2/32.$$
Let $j_1$ and $j_2$ be in $\mathbb{N}^*$ such that
\begin{equation}\label{def-j}
\begin{array}{l}
(c_0/2) n^{ \beta _2/(\beta _1 + \beta _2 + 2\beta _1 \beta _2)} \leq 2^{j_1} \leq c_0  n^{ \beta _2/(\beta _1 + \beta _2 + 2\beta _1 \beta _2)}\\
(c_0/2) n^{ \beta _1/(\beta _1 + \beta _2 + 2\beta _1 \beta _2)} \leq 2^{j_2} \leq c_0  n^{ \beta _1/(\beta _1 + \beta _2 + 2\beta _1 \beta _2)}.
\end{array}
\end{equation}
The existence of $j_1$ and $j_2$ is guaranteed for $n$ larger than an integer $n_0$ depending on $(c_0,\beta)$. 
Then for every $i, i' \in \{0, \dots , M\}$

\begin{eqnarray}\label{eq-distance}
\nonumber \| F_{\varepsilon^{(i)}} - F_{\varepsilon ^{(i')}} \| ^2 &\geq  &\frac{b}{n} \frac{2^{j_1+j_2}}{8}  \geq  \frac{b c_0^2}{32n}  n^{ (\beta _1+\beta _2)/(\beta _1 + \beta _2 + 2\beta _1 \beta _2) }\\
& = & B_0 n^{-2\beta _1 \beta _2/ (\beta _1 + \beta _2 + 2 \beta _1 \beta _2)}  =  B_0 n^{- \overline{\beta} /(\overline{\beta} + 1)}
\end{eqnarray}
hence condition $(2)$ in Theorem $\ref{th-tsy}$ is satisfied. \\

Moreover 
\begin{equation}\label{eq-Kullback}
\frac{1}{M} \sum _{l=0} ^M K( P_{\varepsilon ^{(l)} } ^{(n)} , P_0 ^{(n)})  \leq   a'\| f_{(X,T)}\| _{\infty}  b 2^{j_1+j_2}  =   \frac{\log 2}{72} 2^{j_1+ j_2} \leq   \frac{\log M}{9}
\end{equation}
hence condition $(3)$ in Theorem $\ref{th-tsy}$ is satisfied with $\alpha =1/9$. \\

Finally condition $(1)$ in Theorem $\ref{th-tsy}$ is satisfied as soon as $(\ref{C1})$ and $(\ref{C2})$ hold. Besides, $\beta_1>0$ and $\beta _2>1$ and by ($\ref{def-j}$) $j_1$ and $j_2$ are increasing with $n$. Therefore $2^{(j_1+j_2)/2} (2^{j_1\beta_1} + 2^{j_2\beta_2} +1)$ increases faster than $2^{j_1/2} 2^{3j_2/2} $ and for $n$ larger than an integer $n_1$ depending on $\psi $ and $L$, (\ref{C1}) holds as soon as (\ref{C2}) holds.
Moreover  $(\ref{C2})$
holds as soon as
\begin{equation*}
\sqrt{b} c_0 n^{- \beta _1 \beta _2 /( \beta _1 + \beta _2 + 2\beta _1 \beta _2)} \left( (c_0^{\beta _1} + c_0^{\beta _2} ) n^{\beta _1 \beta _2 /( \beta _1 + \beta _2 + 2\beta _1 \beta _2)} + 1\right) \leq L \left( 1-\sqrt{\frac{7}{12}} \right) \\
\end{equation*}
which is ensured if
\begin{equation} \label{lower-eq1}
\sqrt{b} c_0 (c_0^{\beta _1} + c_0^{\beta _2}) \leq \frac{L}{2} \left( 1-\sqrt{\frac{7}{12}} \right) \quad \text{and}
\end{equation}
\begin{equation}\label{lower-eq2}
\sqrt{b} c_0 n^{- \beta _1 \beta _2 /( \beta _1 + \beta _2 + 2\beta _1 \beta _2)} \leq \frac{L}{2} \left( 1-\sqrt{\frac{7}{12}} \right).
\end{equation}
On the one hand $(\ref{lower-eq1})$ holds as soon as $$2 c_0^{\beta _1+\beta _2+1} \leq \frac{L}{2\sqrt{b}} \left( 1-\sqrt{\frac{7}{12}} \right)$$
which is guaranteed by  the definition of $c_0$. On the other hand there exists an integer $n_2$ depending on $(\beta, c_0)$ such that $(\ref{lower-eq2})$ is satisfied for every $n\geq n_2$.

Thus for every $n\geq \max (n_0, n_1,n_2)$, conditions (1), (2) and (3) in Theorem \ref{th-tsy} hold with $r = B_0  n^{- \overline{\beta} /(\overline{\beta} + 1)}$, which concludes the proof of Theorem \ref{T-low}.  $\Box$

\subsubsection{Minimax rate of convergence in the absence of covariate}
For sake of simplicity, we keep the same notations. Let $(Y_i,T_i) _{i=1, \dots, n}$ be an i.i.d. sample with $T_i$ and $Y_i$ independent positive random variable. Assume that a sample $(T_i,\delta_i)_{i=1,\dots,n}$ is observed with $\delta _i = \text{1\!I}_{\{Y_i\leq T_i\}}$.  Let $F(y)=P[Y_i\leq y]$ the c.d.f. of $Y_i$. Let $\beta>1$,  $L>0$ and $\mathcal{B}^{\beta}_{2, \infty} ([0,1],L)$ the one-dimensional Besov ball of regularity $\beta $ and radium $L$. Then there exists a constant $c$ such that
$$\inf_{\widehat{F}_n } \sup _{F\in \mathcal{B} ^{\beta}_{2, \infty} ([0,1],L)} \mathbb{E} \left[ r \| \widehat{F}_n - F \| ^2 \right] \geq c.$$
The proof is identical to the proof of Theorem \ref{T-low} except that the dependence with respect to $X$ is omitted. More precisely, let $a = \min (1/3,L/2)$ and 
 $$F_0(u)= a1\!{\mbox I}_{[0, +\infty[ }(u) + a u1\!{\mbox I}_{[0,1]}(u) + (1-a)1\!{\mbox I}_{(1, + \infty[} (u) $$
for every $u\in \mathbb{R}$.  Let 
$$b= \frac{ \log 2}{72 \| f_T \| _{\infty} a'}  , \quad  \quad c_0= \left[ \frac{L}{2\sqrt{b}}\left(1- \sqrt{\frac{7}{12}}\right) \right]^{1/(1+ \beta )}  \quad \text{and} \quad B_0=b c_0/16$$
with $a'= (1/a) + (1/(1-2a))$ and $f_T$ the density of $T_i$. Let $j\in \mathbb{N} ^*$ be such that $ (c_0/2)n^{1/(2\beta+1)} \leq 2^j \leq c_0n^{1/(2\beta+1)}.$
For $s\in \mathbb{N} ^*$, let $\psi _{j,s}(u) = 2^{j/2}  \psi (2^ju-s)$. Let $R_j$ be a subset of $\mathbb{Z}$ such that:

\begin{itemize}
\item $Supp(\psi_{j,s})\subset [0,1]$,   $\forall s\in R_j$.
\item The functions $\left\{ \psi _{j,s}, s\in R_j \right\} $ have disjoint supports.
\item $|R_j|=2^j$.
\end{itemize}
Let $M \in \mathbb{N}^*$ and $(\varepsilon ^{(0)}, \dots, \varepsilon ^{(M)}) \in \{0,1\}^{|R_j|}$ be such that $\log(M) \geq (\log 2/8)2^j$ and
$$\overline{\rho}(\varepsilon ^{(i)} , \varepsilon ^{(i')}) \geq \frac{|R_j|}{8}  , \quad \forall i\neq i' \in \{0, \dots ,M\}.$$
For every $i\in \{0, \dots ,M\}$, let 
$$F_{\varepsilon ^{(i)}}  = F_0 + \sqrt{\frac{b}{n}} \sum _{s\in R_j} \varepsilon ^{(i)}_s \psi_{j,s}.$$
Then:

\begin{itemize}
\item Similarly to (\ref{eq-Kullback}), 
\begin{equation*}
\frac{1}{M} \sum _{l=0} ^M K( P_{\varepsilon ^{(l)} } ^{(n)} , P_0 ^{(n)}) \leq   \frac{\log M}{9}
\end{equation*}

\item Similarly to equation (\ref{eq-distance}), for every $i,i' \in \{0,\dots,M\}$
\begin{equation*}
 \| F_{\varepsilon^{(i)}} - F_{\varepsilon ^{(i')}} \| ^2 \geq   B_0 n^{-2\beta/(2\beta + 1)}
\end{equation*}

\item Similarly to (\ref{C1}) and (\ref{C2}), $F_{\varepsilon ^{(i)}} $ is a c.d.f. and lies in $\mathcal{B} ^{\beta}_{2, \infty} ([0,1],L) $ as soon as 

 \begin{eqnarray*} 
 \text{(a)}&& \sqrt{\frac{b}{n}} 2^{3j/2} \| \psi ' \| _{\infty} \leq \frac{a}{2}\\
 \text{(b)} && \sqrt{\frac{b}{n}} 2^{j/2} (2^{j\beta} +1) \leq L\left(  1-\sqrt{\frac{7}{12}} \right)
 \end{eqnarray*}
\end{itemize}
Since $\beta >1$, $(a)$ holds as soon as $(b)$ is satisfied and $(b)$ is guaranteed by the definition of $c_0$.  Thus, the conditions of Theorem \ref{th-tsy} are fulfilled for $r$ proportional to $n^{-2\beta/(2\beta +1)}$. $\Box$

\section*{Acknowledgements}

I am grateful to Fabienne Comte for her advices about this work and to Cristina Butucea for her help on the minimax study.

\bibliographystyle{elsarticle-harv}
\bibliography{Biblio}







\end{document}